     	\def\version{\today}

\newif\iffinal
\finalfalse	
\finaltrue 

\documentclass[reqno,twoside,11pt]{amsart}
\iffinal\else\usepackage[notref,notcite]{showkeys}\fi
\iffinal\else\IfFileExists{pdfsync.sty}{\usepackage{pdfsync}}{}\fi
\usepackage{cite}
\usepackage{amsmath}
\usepackage{amsfonts}
\usepackage{amssymb}
\usepackage{epsfig}
\usepackage{verbatim}
\usepackage{fancyhdr}
\usepackage{mathabx}

\IfFileExists{epsf.def}{\input epsf.def}{\usepackage{epsf}}

\usepackage{mathrsfs}
\usepackage{mathpazo}
\usepackage{mathabx}

\DeclareFontFamily{OT1}{eusb}{} \DeclareFontShape{OT1}{eusb}{m}{n} {<5> <6> <7> <8> <9> <10> <11> <12> <14.4> eusb10}{}
\DeclareMathAlphabet{\eusb}{OT1}{eusb}{m}{n}

\DeclareFontFamily{OT1}{eusm}{} \DeclareFontShape{OT1}{eusm}{m}{n} {<5> <6> <7> <8> <9> <10> <11> <12> <14.4> eusm10}{}
\DeclareMathAlphabet{\eusm}{OT1}{eusm}{m}{n}

\DeclareFontFamily{OT1}{eufm}{} \DeclareFontShape{OT1}{eufm}{m}{n} {<5> <6> <7> <8> <9> <10> <11> <12> <14.4> eufm10}{}
\DeclareMathAlphabet{\mathfrak}{OT1}{eufm}{m}{n}

\DeclareFontFamily{OT1}{fraktura}{}
\DeclareFontShape{OT1}{fraktura}{m}{n} {<5> <6> <7> <8> <9> <10> <11> <12> <13> <14.4> [1.1] eufm10}{}
\DeclareMathAlphabet{\fraktura}{OT1}{fraktura}{m}{n}

\DeclareFontFamily{OT1}{cmfi}{} \DeclareFontShape{OT1}{cmfi}{m}{n} {<5> <6> <7> <8> <9> <10> <11> <12> <13> <14.4> [0.9] cmfi10}{}
\DeclareMathAlphabet{\cmfi}{OT1}{cmfi}{b}{n}

\DeclareFontFamily{OT1}{cmss}{} \DeclareFontShape{OT1}{cmss}{m}{n} {<5> <6> <7> <8> <9> <10> <11> <12> <13> <14.4> cmss10}{}
\DeclareMathAlphabet{\cmss}{OT1}{cmss}{m}{n}

\newcommand{\lefttinyfoottext}{}
\newcommand{\righttinyfoottext}{\hfill Preliminary version (typeset: \version)}

\fancypagestyle{plain}{%
  \fancyhf{}%
  \fancyfoot[LE,LO]{{\tiny \lefttinyfoottext}}
\fancyfoot[RE,RO]{{\tiny \righttinyfoottext}}
\fancyfoot[CE,CO]{\tiny\thepage}

}

\newtheoremstyle{thm}{1.5ex}{1.5ex}{\itshape\rmfamily}{} {\bfseries\rmfamily}{}{2ex}{}

\newtheoremstyle{def}{1.5ex}{1.5ex}{\rmfamily\sl}{} {\bfseries\rmfamily}{}{2ex}{}

\newtheoremstyle{rem}{1.5ex}{1.5ex}{\rmfamily}{} {\bfseries\rmfamily}{}{2ex}{}

\newenvironment{proofsect}[1] {\vskip0.1cm\noindent{\bfseries\itshape#1.}}{\qed\vspace{0.15cm}}

\theoremstyle{thm}
\newtheorem{theorem}{Theorem}[section]
\newtheorem{lemma}[theorem]{Lemma}

\newtheorem*{Main Theorem}{Main Theorem.}
\newtheorem{corollary}[theorem]{Corollary}
\newtheorem{exercise}[theorem]{Exercise}
\newtheorem{conjecture}[theorem]{Conjecture}

\theoremstyle{def}
\newtheorem{definition}[theorem]{Definition}

\theoremstyle{rem}
\newtheorem{remark}[theorem]{{Remark}}

\newtheorem{question}[theorem]{Question}

\numberwithin{equation}{section}

\newcommand{\sectionphrase}{Lecture~\arabic{section}\ }

\renewcommand{\section}{\secdef\sct\sect}
\newcommand{\sct}[2][default]{\refstepcounter{section}
\addcontentsline{toc}{section}
{{\tocsection {}{\thesection}{\!\!\!\!#1\dotfill}}{}}
\vspace{0.7cm}
\vglue2cm
\leftline{\bfseries\Large\sectionphrase}
\vglue6mm
\leftline{\bfseries\LARGE #1} \nopagebreak \vglue1cm}
\newcommand{\sect}[1]{
\vspace{0.4cm} \centerline{\large\scshape\rmfamily #1}
\vspace{0.2cm}}

\renewcommand{\subsection}{\secdef\subsct\sbsect}
\newcommand{\subsct}[2][default]{\refstepcounter{subsection}
\addcontentsline{toc}{subsection}
{{\tocsection{\!\!}{\hspace{1.2em}\thesubsection}{\!\!\!\!#1\dotfill}}{}}
\nopagebreak\vspace{0.45\baselineskip} {\flushleft\bf
\arabic{section}.\arabic{subsection}~\bf #1.~}
\\*[3mm]\noindent
\nopagebreak}
\newcommand{\sbsect}[1]{\vspace{0.1cm}\noindent
\textbf{#1.~}\vspace{0.1cm}}

\renewcommand{\subsubsection}{%
\secdef \subsubsect\sbsbsect}
\newcommand{\subsubsect}[2][default]{%
\refstepcounter{subsubsection} 
\addcontentsline{toc}{subsubsection}{{\tocsection{\!\!}
{\hspace{3.05em}\thesubsubsection}{\!\!\!\!#1\dotfill}}{}}
\nopagebreak
\vspace{0.15\baselineskip} \nopagebreak {\flushleft\rmfamily
\itshape\arabic{section}.\arabic{subsection}.\arabic{subsubsection}
\ \rmfamily #1\/.}\ }
\newcommand{\sbsbsect}[1]{\vspace{0.1cm}\noindent
\rmfamily \itshape
\arabic{section}.\arabic{subsection}.\arabic{subsubsection} \
\sffamily #1\/.\ }

\iffinal

\else

\fi

\renewcommand{\caption}[1]{%
\vglue0.5cm
\refstepcounter{figure}
\begin{minipage}{0.9\textwidth}\small {\sc Figure~\thefigure. }#1\end{minipage}}

\newcommand{\supp}{\operatorname{supp}}

\newcommand{\textd}{\text{\rm d}\mkern0.3mu}

\newcommand{\texte}{\text{\rm  e}\mkern0.7mu}

\newcommand{\Var}{\text{\rm Var}}
\newcommand{\Cov}{\text{\rm Cov}}

\renewcommand{\AA}{\mathcal A}

\newcommand{\DD}{\mathcal D}

\newcommand{\NN}{\mathcal N}

\newcommand{\TT}{\mathcal T}

\newcommand{\E}{\mathbb E}

\newcommand{\N}{\mathbb N}

\newcommand{\BbbP}{\mathbb P}

\newcommand{\R}{\mathbb R}

\newcommand{\Z}{\mathbb Z}

\newcommand{\twoeqref}[2]{(\ref{#1}--\ref{#2})}

\def\myffrac#1#2 in #3{\raise 2.6pt\hbox{$#3 #1$}\mkern-1.5mu\raise 0.8pt\hbox{$#3/$}\mkern-1.1mu\lower 1.5pt\hbox{$#3 #2$}}

\newcommand\independent{\protect\mathpalette{\protect\independenT}{\perp}}
\def\independenT#1#2{\mathrel{\rlap{$#1#2$}\mkern3mu{#1#2}}}

\newcommand{\wh}{\widehat}
\newcommand{\wt}{\widetilde}
\newcommand{\cc}{{\text{\rm c}}}

\newcommand{\laweq}{\overset{\text{\rm law}}=}

\newcommand{\lawarrow}{{\overset{\text{\rm law}}\longrightarrow}}
\newcommand{\fraka}{\fraktura a}
\newcommand{\frakg}{\fraktura g}
\newcommand{\frakd}{\fraktura d}

\begin{document}

\vglue-5mm
\title[\theruntitle]{{\rm Lecture notes for the minicourse}\\*[4mm]\large Extremal properties of the random walk local time\\*[4mm]{\rm \normalsize School on disordered media\\*[2mm]Alfr\'ed R\'enyi Insititute, Budapest\\
January 20-24, 2025}}
\author[]{\bf\large Marek Biskup$^\ast$}

\maketitle
\thispagestyle{empty}

\vglue-3mm
\vbox{
\setcounter{tocdepth}{2}
\tableofcontents
}

\vglue-3mm

\vfill

\leftline{\small$^\ast$Address: Department of Mathematics, UCLA, Los Angeles, CA 90095-1555, U.S.A.}
\leftline{\small\copyright 2025 M.~Biskup. Reproduction for non-commercial use is permitted.}
\leftline{\small Typeset: \today}




\newpage

\setcounter{page}{1}

\section{Introduction and main results}
\noindent
The minicourse to be given over four 50-minute lectures will focus on extremal properties of random walk local time. This turns out to be a particular aspect of the larger area of logarithmically correlated processes that has attracted a lot of attention in recent years. For lack of time, we will focus only on one particular result; namely, the scaling limit of the points avoided by a two-dimensional simple random walk. The main objective of the course is to motivate the students to learn other, and often more difficult, results through self-study of papers and existing review articles.

\subsection{Random walk local time}
Throughout we will consider a continuous time Markov chain~$X$ on a finite state space that in general will take the form $V\cup\{\varrho\}$, where $V$ is a finite set and~$\varrho$ is a distinguished vertex (not belonging to~$V$). The transitions will occur at ``constant speed,'' which means that the chain takes the form $X_t = Z_{N(t)}$, where~$Z$ is a discrete-time Markov chain and~$\{N(t)\colon t\ge0\}$ is a rate-$1$ Poisson point process independent of~$Z$. We assume that~$X$ (and~$Z$) is irreducible and reversible with respect to measure~$\pi$ and write~$P^x$ for the law of the chain started at~$x$, with associated expectation denoted as~$E^x$.

Of our prime interest in these lectures is the \textit{local time} associated with~$X$. This is the two-parameter stochastic  process
\begin{equation}
\bigl\{\ell_t(x)\colon x\in V\cup\{\varrho\},\,t\ge0\bigr\}
\end{equation}
 defined by
\begin{equation}
\ell_t(x):=\frac1{\pi(x)}\int_0^t 1_{\{X_s=x\}}\textd s
\end{equation}
(Recall that irreducibility forces $\pi(x)>0$ for all~$x\in V\cup\{\varrho\}$.) To explain the normalization, note that $X_t$ will for large~$t$ be distributed according to suitably-normalized~$\pi$ which means that the time spent at~$x$ grows proportionally to~$\pi(x)$.

The main question of interest in these lectures is then:
\begin{equation}
\text{What does $\ell_t$ look like at large~$t$?}
\end{equation}
While this is our general objective, we focus on particular questions. For example, we may ask about the size and asymptotic law of
\begin{equation}
\max_{x\in V\cup\{\varrho\}}\ell_t(x)\,\,\text{ and }\,\min_{x\in V\cup\{\varrho\}}\ell_t(x)
\end{equation}
Noting that the minimum is zero until the first time all vertices are visited naturally leads us to the notion of the \textit{cover time},
\begin{equation}
\tau_{\text{cov}}:=\inf\Bigl\{t\ge0\colon \min_{x\in V\cup\{\varrho\}}\ell_t(x)>0\Bigr\}
\end{equation}
of which we can then ask how it scales with~$t$ and the size of~$V$.
Another natural question (which is the one we will focus in these notes) concerns the structure of the set of points not yet visited by~$X$ at time~$t$; namely,
\begin{equation}
\label{E:1.5}
\AA(t):=\bigl\{x\in V\cup\{\varrho\}\colon \ell_t(x)=0\bigr\}
\end{equation}
that we will refer to as \textit{avoided points}.

Of course, taking~$t$ to be large without changing~$V$ will hardly lead to interesting conclusions. We will also treat only one particular class of Markov chains; namely, that arising from the simple symmetric random walk (SRW) on~$\Z^d$. So, unless we discuss general aspects of the theory where the setting introduced above is more appropriate, we take~$V$ to be a scaled-up and discretized version~$D_N\subseteq\Z^d$ of a ``nice'' continuum domain $D\subseteq\R^d$; i.e., the set (roughly) of the form
\begin{equation}
\label{E:1.6}
\bigl\{x\in\Z^d\colon x/N\in D\bigr\}
\end{equation}
or (in some references to prior work) the lattice torus~$(\Z/N\Z)^d$. The point is now to study the local time for the simple random walk in~$D_N$ at times~$t_N$ such that~$t_N\to\infty$, subject to specific restrictions on growth rate with~$N$ as~$N\to\infty$.

A technical point for~$D_N$ of the form \eqref{E:1.6} is how to interpret the ``random walk'' at the vertices of~$D_N$ that, in~$\Z^d$, would have an edge to the complement of~$D_N$. One possibility is to treat this as a \textit{free boundary condition} which means to ignore jumps leading out of~$D_N$. For reasons that will become clear later we use a different ``return mechanism,'' corresponding to the \textit{wired boundary condition}, which is defined as follows: Collapse all the vertices in~$\Z^d\smallsetminus D_N$ to one \textit{boundary vertex}~$\varrho$. Then route the edges emanating out of~$D_N$ to~$\varrho$. This leads to a domain as in Fig.~\ref{fig1}.

\medskip
\newcounter{obrazek}
\refstepcounter{obrazek}
\centerline{\hglue1cm$D_N$}
\centerline{\includegraphics[width=2.8in]{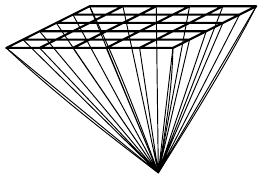}}
\vglue-5mm
\centerline{\hglue1.9cm$\varrho$}
\begin{quote}
\label{fig1}
\fontsize{9}{5}\selectfont
{\bf Fig.~\theobrazek:\ }An illustration of the state space for the random walk. Here~$D_N$ is simply an $N\times N$ square while~$\varrho$ is a vertex to which all the boundary edges of~$D_N$ in~$\Z^2$ are re-routed.
\normalsize
\medskip
\end{quote}

\begin{figure}[t]
\refstepcounter{obrazek}
\centerline{\includegraphics[width=3.4in]{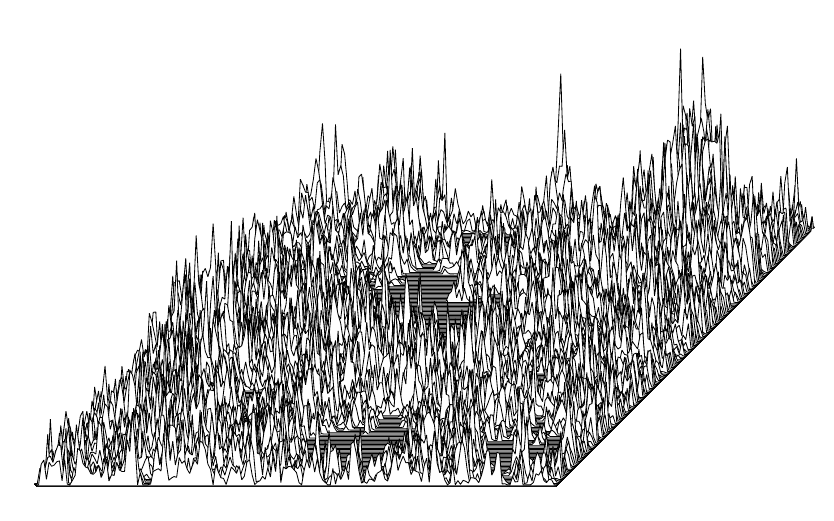}
\hglue4mm \includegraphics[width=1.8in]{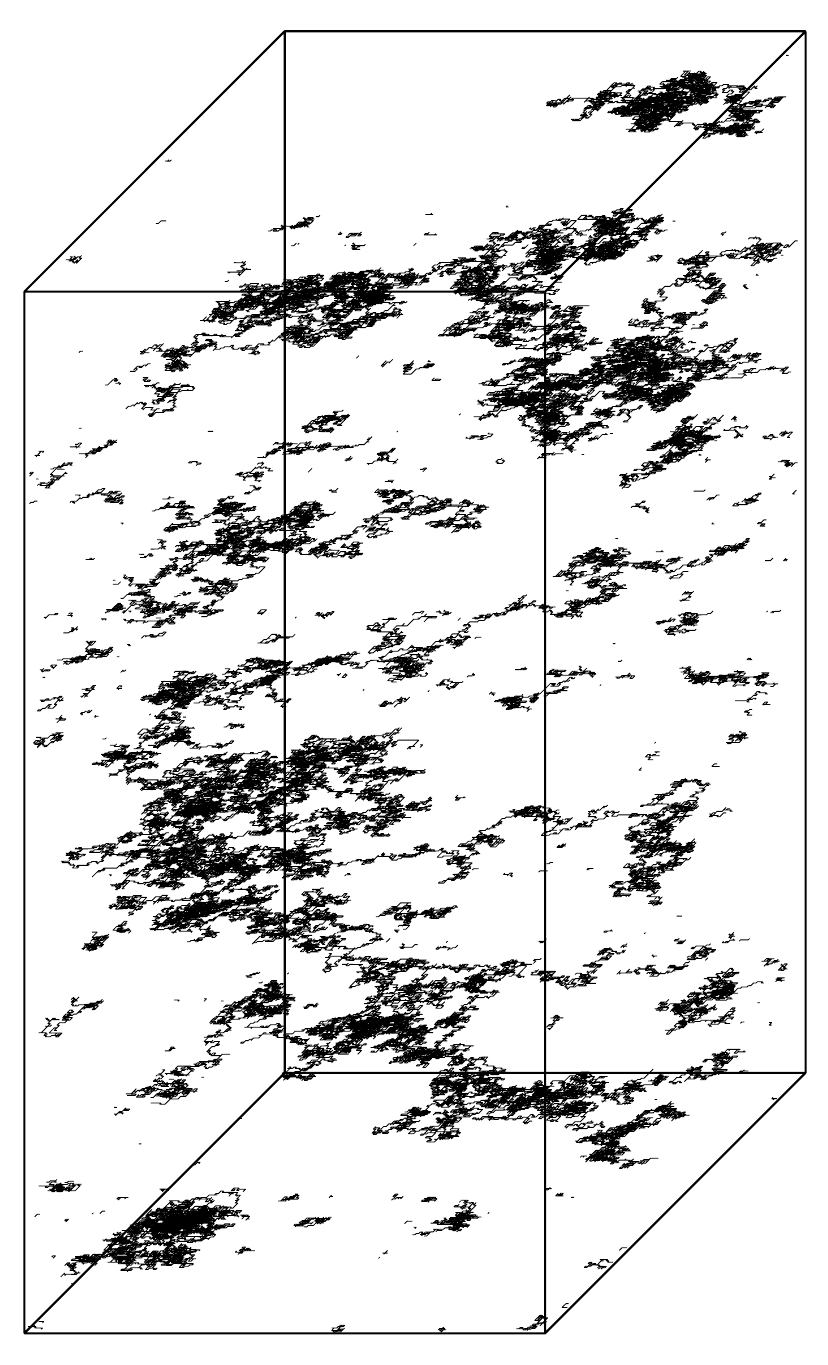}}
\smallskip
\begin{quote}
\label{fig2a}
\fontsize{9}{5}\selectfont
{\bf Fig.~\theobrazek:\ }A sample of the local time (left) and the trajectory of the walk (right) over time (with time axis running upwards) for the random walk on~$D_N\cup\{\varrho\}$ as in~Fig.~\ref{fig1} with~$N=200$.
\normalsize
\medskip
\end{quote}
\end{figure}

\noindent
Note that the invariant measure~$\pi(x)$ equals the degree of the vertex~$x$ in the resulting graph which for~$D_N$ with wired boundary condition will be equal~$2d$ at $x\in D_N$ and equal to the size of the edge boundary of~$D_N$ in~$\Z^d$ at~$\varrho$. The chain~$X$ is then a constant-speed continuous-time random walk on the resulting finite graph that runs just as the simple random walk on~$D_N$ and, after each exit, returns back to~$D_N$ through a randomly-chosen boundary edge. (For readers worried that this might lead to the local time building up near the boundary, Fig.~\ref{fig2a} and our Theorem~\ref{thm-1.7} show that this is not the case.)

\subsection{The case for~$d=2$}
In these lectures we will focus on the above setting in spatial dimension~$d=2$. To motivate this, let us recount some of the basic facts about random walk on~$\Z^d$. The key difference arises already in the celebrated P\'olya theorem that says
\begin{equation}
\text{SRW on~$\Z^d$ is }\begin{cases}
\text{transient},\qquad&\text{in }d\ge3
\\
\text{recurrent},\qquad&\text{in }d=1,2
\end{cases}
\end{equation}
The transience can be thought of as a ``short memory'' (or decaying-autocorrelation) property that very often makes a number of arguments easier to handle. 

The discrepancy between the recurrent and transient regime typically manifests itself in the analytic form of the conclusions. To give an example, note that for the cover time of the lattice torus~$(\Z/N\Z)^d$ the following holds:
\begin{equation}
\label{E:1.8}
\tau_{\text{cov}}\asymp\begin{cases}
N^d\log N,\qquad&\text{in }d\ge3
\\
N^2(\log N)^2,\qquad&\text{in }d=2
\\
N^2,\qquad&\text{in }d=1
\end{cases}
\end{equation}
where the first two lines are true as sharp asymptotics (with a known constant of proportionality) because the cover time concentrates strongly around its expected value. This concentration fails in~$d=1$ as $\tau_{\text{cov}}/N^2$ tends in law (as~$N\to\infty$) to a non-degenerate random variable. 

The formula~$N^d\log N$ in~$d\ge3$ is easy to understand: One needs~$N^d$ time to visit most of the vertices but then a coupon-collector reasoning need to be employed to sweep out the  outliers. A similar reasoning can be used in~$d=2$; the extra~$\log N$ appears because once a vertex is hit, it will be visited order~$\log N$ times before it is left for good. (Much more is known in fact; thanks to a result of D.~Belius~\cite{Belius}, we know a full limit law of suitably centered $\tau_{\text{cov}}$ in all~$d\ge3$. In~$d=2$, the corresponding asymptotic is the subject of active research by several groups.)

Moving to the set of avoided points, the natural time scales to consider are those proportional to the cover time. So we will specialize \eqref{E:1.5} to the form
\begin{equation}
\AA_N(\theta):=\bigl\{x\in D_N\colon\ell_{\theta E^\varrho\tau_{\text{cov}}}(x)=0\bigr\}
\end{equation}
We now ask about the asymptotic properties of this set, specifically, the number of points and the way they are distributed in~$D_N$, in the limit as~$N\to\infty$. 

Also in this problem the recurrence/transience dichotomy manifests itself strong\-ly in the conclusions. Indeed, in $d\ge3$ the set $\AA_N(\theta)$ partitions into a collection of independent small ``islands'' which, thanks to a result of J.~Miller and P.~Sousi~\cite{Miller-Souzi} from 2017 can even be nailed to the form
\begin{equation}
\AA_N(\theta) \,\,\overset{\text{a.s.}}\approx\,\, \text{Bernoulli}(N^{-\theta d})
\end{equation}
for a non-trivial interval of~$\theta\in[0,1]$. Here the squiggly equality represents a coupling in total variational distance to the set where the Bernoulli process equals~$1$. 

\begin{figure}[t]
\refstepcounter{obrazek}
\centerline{\includegraphics[width=2.2in]{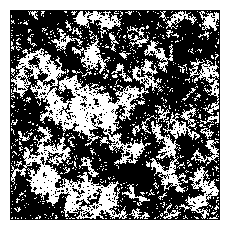}
\hglue4mm \includegraphics[width=2.2in]{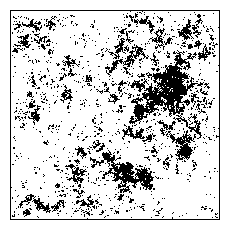}}
\begin{quote}
\label{fig2}
\fontsize{9}{5}\selectfont
{\bf Fig.~\theobrazek:\ }Samples of the set of avoided points for the random walk on~$D_N$ run for time proportional to~$\theta=0.1$ (left) and~$\theta=0.3$ (right) fraction of the expected cover time. (The same run of the random walk is used for both figures.)
\normalsize
\medskip
\end{quote}
\end{figure}

In contrast to this, in $d=2$, the set~$\AA_N(\theta)$ scales to a \textit{random fractal}, as shown in Fig.~\ref{fig2}. The point of these lectures is to make sense of a limit of these pictures as~$N\to\infty$.

\subsection{Link to Gaussian Free Field}
Our method to control the local time will rely on a close connection between the local time of a Markov chain and a Gaussian process called \textit{Gaussian Free Field}. We will now introduce this concept in the general setting of Markov chains on~$V\cup\{\varrho\}$ introduced above. The connection itself will be discussed in Lecture~2.

We start with some standard definitions. For any~$x\in V\cup\{\varrho\}$ we introduce the \textit{first hitting time} of~$x$ by~$X$ as
\begin{equation}
H_x:=\inf\{t\ge0\colon X_t=x\}
\end{equation}
Notice that~$H_x=0$ $P^x$-a.s. We then use this to define the \textit{Green function}
\begin{equation}
G^V\colon (V\cup\{\varrho\})\times(V\cup\{\varrho\})\to[0,\infty)
\end{equation}
by the formula
\begin{equation}
G^V(x,y):= E^x\bigl(\ell_{H_\varrho}(y)\bigr) = \frac1{\pi(y)}\int_0^\infty P^x\bigl(X_t=y,\,H_\varrho> t\bigr)\textd t
\end{equation}
where the second expression is based on writing $\ell_{H_\varrho}(y) = \frac1{\pi(y)}\int_0^{H_\varrho} 1_{\{X_t=t\}}\textd t$ and applying Tonelli's theorem. We now pose our first exercise:

\begin{exercise}
Show that, if viewed as a matrix, $G^V$ is symmetric and positive semidefinite.
\end{exercise}

\noindent
The reason why we highlight these properties is that they make~$G^V$ a covariance. This is enough to make sense of:

\begin{definition}
The Discrete Gaussian Free Field (DGFF) on~$V$ is a Gaussian process
\begin{equation}
\{h^V_x\colon x\in V\cup\{\varrho\}\}
\end{equation}
 such that
\begin{equation}
\forall x,y\in V\cup\{\varrho\}\colon\quad \E h^V_x=0\quad\wedge\quad \E(h^V_xh^V_y)=G^V(x,y)
\end{equation}
(We will use~$\BbbP$ and~$\E$ for probability and expectation associated with the~DGFF.)
\end{definition}

Note that the definitions ensure that $G^V(\varrho,y)=0=G^V(x,\varrho)$ for any~$x$ and~$y$. This along with $\E h_\varrho=0$ forces
\begin{equation}
h^V_\varrho = 0\quad\BbbP\text{-a.s.}
\end{equation}
which also explains the special role the ``boundary vertex''~$\varrho$ plays in the whole setup. In particular, our~$h^V$ corresponds to the case of Dirichlet boundary conditions.

The reason for calling this the ``Discrete'' GFF is to make a distinction between the corresponding concept in the continuum, called the ``Continuum'' GFF with the shorthand~CGFF. While the latter is not a prime target of interest in our notes, we will make some references to it when we discuss the scaling limits in Lectures~3 and~4.

When we specialize ourselves to the random walk on $V:=D_N$ and take the resulting DGFF $h^{D_N}$ at face value, we are naturally led to ask a number of questions about its extremal properties similar as those for the local time asked above. For instance:
\begin{equation}
\text{What is the growth rate/scaling limit of }\max_{x\in D_N}h^{D_N}_x?
\end{equation}
(Since the field is symmetric, the minimum scales as the negative of the maximum.) Another question to ask is:
\begin{equation}
\text{What is the cardinality/scaling limit of }\Bigl\{x\in D_N\colon h^{D_N}_x\ge\lambda\max_{x\in D_N}h^{D_N}_x\Bigr\}?
\end{equation}
Here, for~$\lambda\in(0,1)$, we call the points in the set the \textit{$\lambda$-thick points} of~$h^{D_N}$.

\begin{figure}[t]
\refstepcounter{obrazek}
\centerline{\includegraphics[width=3.4in]{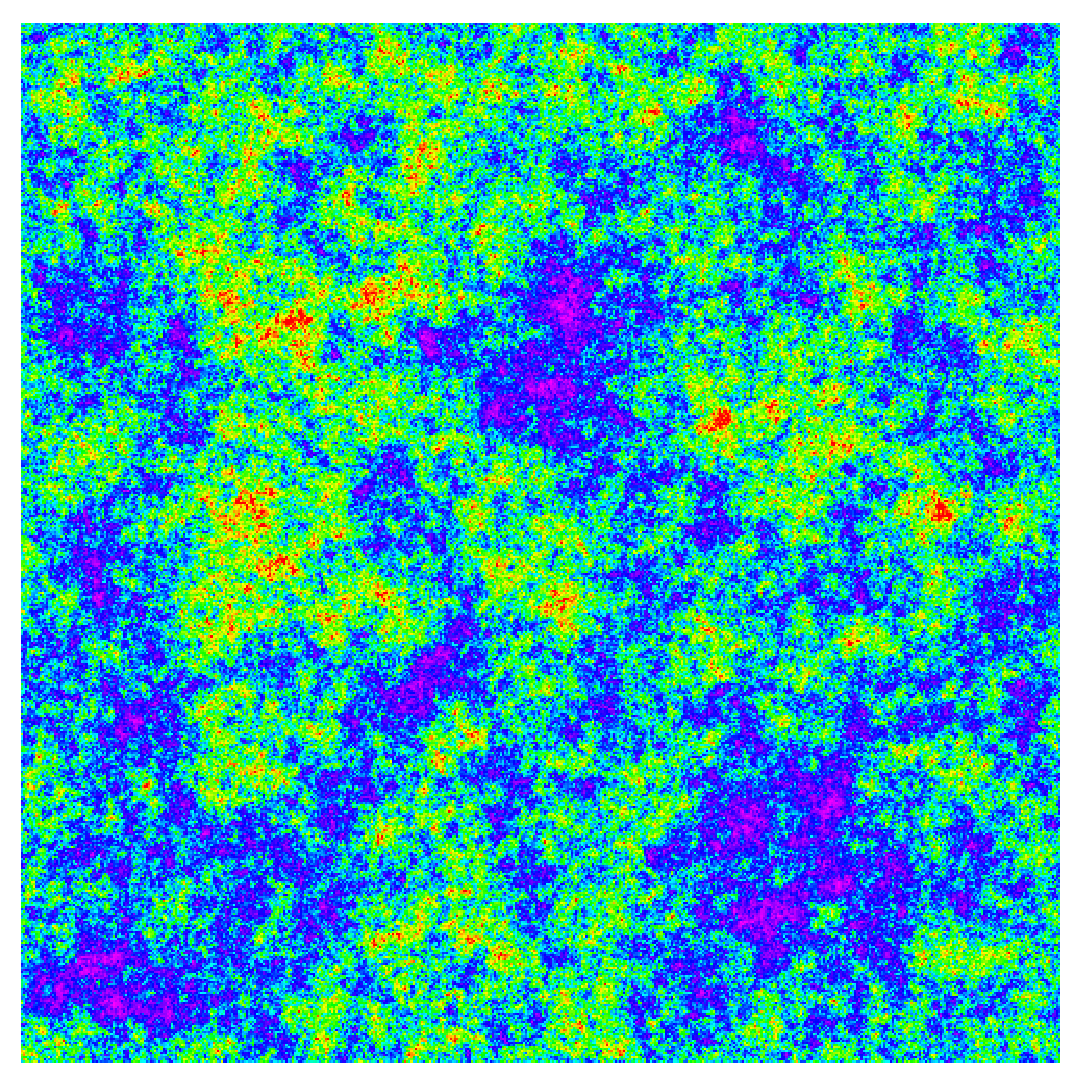}}
\begin{quote}
\label{fig3}
\fontsize{9}{5}\selectfont
{\bf Fig.~\theobrazek:\ }A sample of DGFF on $500\times500$-square color coded so that the red regions are those with large positive values and purple regions are those with large negative values. The values in-between are coded according to usual ordering of colors by the wave-length.
\normalsize
\end{quote}
\vglue-1mm
\end{figure}

Similarly as for the local time, the DGFF samples in spatial dimensions $d\ge3$ are not nearly so interesting as in spatial dimension two. Indeed, in~$d=2$ the field itself is a random fractal which we demonstrate in Fig.~\ref{fig3}. 

Looking at the figure, the reader will surely notice the fractal curves separating the mostly-green and mostly-blue regions; these are known to be described by the SLE$_4$-curves thanks to a celebrated work by O.~Schramm and S.~Sheffield~\cite{Schramm-Sheffield}. Our interest in the present text are the yellow-to-red regions, where the field is unusually large.

A key problem in making any of the above mathematically reasonable is the fact that the DGFF in $d=2$ becomes increasingly singular as the side of the underlying domain increases. This has to do with the fact that the field is \textit{logarithmically correlated}; we will elaborate on what this means in Lecture~2.

\subsection{Main result on DGFF}
We are now ready to make precise statements of the main results to be discussed in detail throughout the rest of the course. We start with those for the DGFF. First we identify the continuum regions~$D$ to which our results apply:

\begin{definition}[Admissible domain]
A set $D\subseteq\R^2$ is an admissible domain if it is bounded, open and $\partial D$ has a finite number of connected components each of which has a positive Euclidean diameter.
\end{definition}

Since a bounded simply connected open subset of~$\R^2$ has a connected boundary, it follows that any such set is an admissible domain by the above definition. However, we also allow for non-trivial arcs in the interior. While connectedness is not required, the fact that our processes will trivially factor over connected components means that we only need to work with connected admissible~$D$.

Next we will specify more precisely the way we allow ourselves to discretize~$D$. While \eqref{E:1.6} seems to be a canonical choice, the problem is that this choice may result in a discrete set that ``looks'' quite different than~$D$ itself; particularly, from the perspective of harmonic analysis. Writing $d_\infty$ for the infinity distance on~$\Z^2$, we instead use:

\begin{definition}[Admissible approximations]
A sequence $\{D_N\}_{N\ge1}$ of non-empty subsets of~$\Z^2$ is an admissible approximation of an admissible domain~$D\subseteq\R^d$ if 
\begin{equation}
\forall N\ge1\colon\quad D_N\subseteq\bigl\{x\in\Z^d\colon d_\infty(x/N,\R^2\smallsetminus D)>1/N\bigr\}
\end{equation}
and, for all $\delta>0$ there exists $N_0\ge1$ such that
\begin{equation}
\forall N\ge N_0\colon\quad D_N\supseteq\bigl\{x\in\Z^d\colon d_\infty(x/N,\R^2\smallsetminus D)>\delta\bigr\}
\end{equation}
\end{definition}

\begin{figure}[t]
\refstepcounter{obrazek}
\centerline{\includegraphics[width=2.8in]{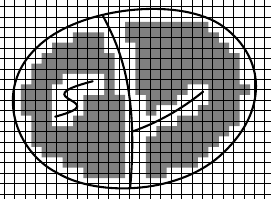}}
\bigskip
\begin{quote}
\label{fig4}
\fontsize{9}{5}\selectfont
{\bf Fig.~\theobrazek:\ }An illustration of an admissible approximation~$D_N$ (marked by the lattice side in the dark region) of an admissible domain~$D\subseteq\R^2$ (bounded by the thick lines).
\normalsize
\end{quote}
\vglue0mm
\end{figure}

\noindent
To illustrate this on an example, Fig.~\ref{fig4} shows an admissible lattice approximation of an admissible domain.
In Lecture~2 we will give the reasons why definitions need to be set up this way. 

In order to set the scales, next we note that for~$x$ deep inside~$D_N$ we will have
\begin{equation}
G^{D_N}(x,x)=g\log N+O(1)
\end{equation}
where the constant of proportionality equals
\begin{equation}
\label{E:1.21}
g:=\frac1{2\pi}
\end{equation}
For the maximum of the DGFF we then get the asymptotic
\begin{equation}
\max_{x\in D_N}h^{D_N}_x=2\sqrt g\,\log N+O(\log\log N)
\end{equation}
Interestingly, the same asymptotic applies even for i.i.d.\ normals with variance $g\log N$; however, while the leading order is the same in the two cases, the constant multiplying $\log\log N$ is already different.

\begin{figure}[t]
\bigskip
\refstepcounter{obrazek}
\centerline{\includegraphics[width=2.2in]{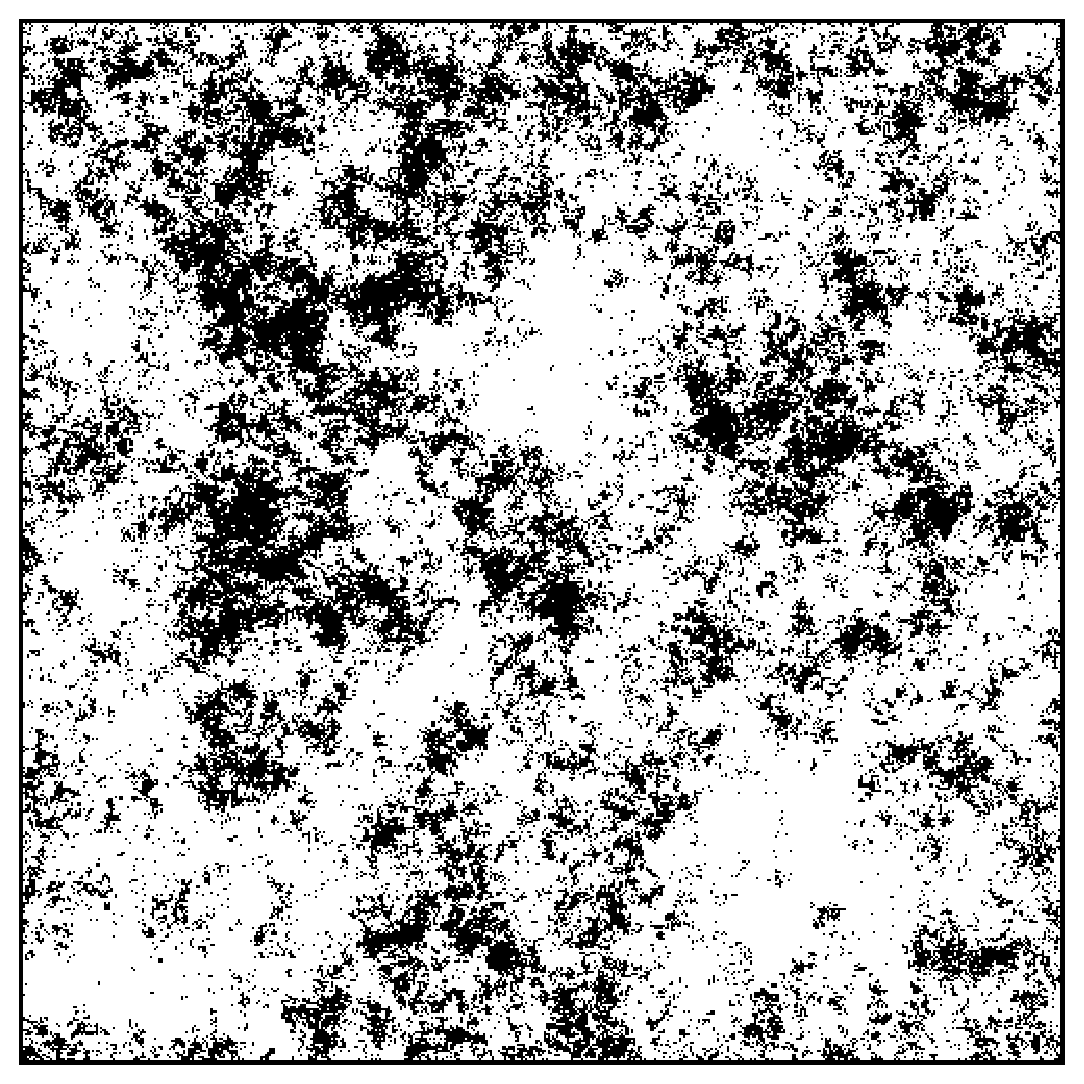}
\hglue4mm \includegraphics[width=2.2in]{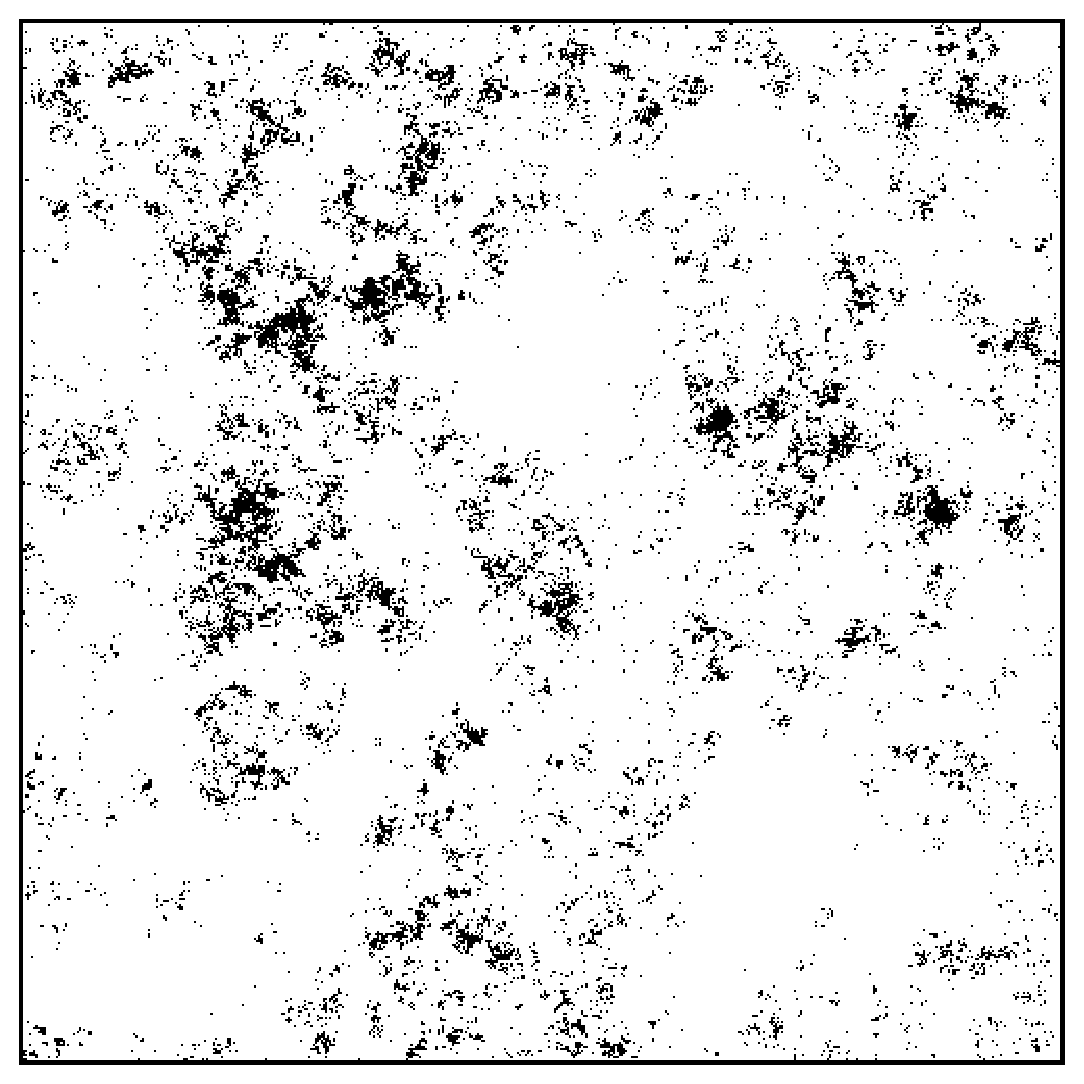}}
\begin{quote}
\label{fig5}
\fontsize{9}{5}\selectfont
{\bf Fig.~\theobrazek:\ } The level sets at $\lambda:=0.1$ (left) and $\lambda=0.3$ (right) multiple of the expected maximum of DGFF on a square box of side-length $500$..
\normalsize
\end{quote}
\vglue4mm
\end{figure}

We now define the set of $\lambda$-thick points again as
\begin{equation}
\label{E:1.23}
\TT_N(\lambda):=\bigl\{x\in D_N\colon h^{D_N}_x\ge  2\sqrt \lambda g\,\log N\bigr\},\quad \lambda\in(0,1).
\end{equation}
As noted earlier, this set is expected to look like a random fractal which readily confirmed by simulation, see Fig.~\ref{fig5}.

A question to address next is in what sense  we can take a scaling limit of the pictures in Fig.~\ref{fig5}. The natural choice of Hausdorff distance is out because, after scaling by~$1/N$ and taking $N\to\infty$, these sets are everywhere dense and so have vanishing asymptotic distance to domain~$D$ itself. We will therefore use a different approach: We associate with each set a point process that records both the position and the value of the field and then take the limit of this process itself.

For instance, for the set \eqref{E:1.23} the point process would take the form
\begin{equation}
\label{E:1.26u}
\sum_{x\in D_N} \delta_{x/N}\otimes\delta_{h^{D_N}_x-2\sqrt g\,\lambda\log N}
\end{equation}
where the tensor product of delta-measures is just a convenient way to write a delta-measure at the corresponding two-coordinate quantity. The point is that, while having total mass of~$D_N$, this measure  does give us access to the cardinality of~$\TT_N(\lambda)$ by integrating it against the function $f(x,h):=1_{[0,\infty)}(h)$. Of course, \eqref{E:1.26u} by itself would not allow for a reasonable limit as $N\to\infty$ as one still needs to normalize the measure in such a way that a limit in law is possible. This is what we do in:

\begin{theorem}[B.-Louidor~\cite{BL4}]
\label{thm-1.5}
Let~$\{D_N\}_{N\ge1}$ be admissible approximations of an admissible domain~$D\subseteq\R^2$. There exists a family of a.s-finite random Borel measures $\{Z^D_\lambda\colon \lambda\in(0,1)\}$ on~$D$ such that for any positive sequence $\{a_N\}_{n\ge1}$ with
\begin{equation}
\label{E:1.25}
\lambda:=\lim_{N\to\infty}\frac{a_N}{2\sqrt g\log N}\in(0,1)
\end{equation}
and with
\begin{equation}
\label{E:1.26}
K_N:=\frac{N^2}{\sqrt{\log N}}\texte^{-\frac{a_N^2}{2g\log N}}
\end{equation}
we have
\begin{equation}
\label{E:1.29u}
\frac1{K_N}\sum_{x\in D_N} \delta_{x/N}\otimes\delta_{h^{D_N}_x-a_N}\,\,\,\underset{N\to\infty}\lawarrow\,\,\,\,Z^D_\lambda(\textd x)\otimes \texte^{-\alpha\lambda h}\textd h
\end{equation}
where $\alpha:=2/\sqrt g$ and the weak convergence is relative to vague topology on $\overline D\times\R$.
Moreover, a.e.\ sample of $Z^D_\lambda$ charges every non-empty open set and is diffuse (i.e., has no atoms).
\end{theorem}

A few remarks are on order: First, the statement \eqref{E:1.29u} means that integrating the measure on the left with respect to any number of continuous functions $\overline D\times\R\to\R$ with compact support results in a family of random variables that converge jointly in law to the corresponding integrals with respect to the measure on the right-hand side. As it turns out, due to linearity of the integral and the Cram\`er-Wold device, it suffices to check convergence for integrals of just one function at the time.

Next observe that \twoeqref{E:1.25}{E:1.26} give
\begin{equation}
K_N = N^{2(1-\lambda^2)+o(1)},\quad N\to\infty.
\end{equation}
Although the total mass of the measure  in \eqref{E:1.25} diverges proportionally to $N^{2\lambda^2+o(1)}$, the mode of convergence ensures that the mass it puts on any compact subset of $\overline D\times\R$ is tight and admits a distributional limit. 

Third, as a corollary of Theorem~\ref{thm-1.5} we get a limit law for the total size of the level set of~$\lambda$-thick points:

\begin{corollary}
For the setting of Theorem~\ref{thm-1.5},
\begin{equation}
\frac1{K_N}\#\bigl\{x\in D_N\colon h^{D_N}_x\ge a_N\bigr\}\,\,\,\underset{N\to\infty}\lawarrow\,\,\,\,(\alpha\lambda)^{-1}\,Z^D_\lambda(D)
\end{equation}
\end{corollary}

\noindent
This generalizes a result of O.~Daviaud~\cite{Daviaud} from 2006 who obtained the leading order growth of the level set without identifying the subleading terms and/or a limit law. The above results will be discussed and, for $\lambda<1/\sqrt2$, proved in Lecture~3.

\subsection{Main result on local time}
Moving to our results on the local time, consider the random walk on~$D_N\cup\{\varrho\}$ as described above. We will for simplicity focus only on the avoided points at times proportional to the cover time of~$D_N\cup\{\varrho\}$. As it turns out, the easiest format to state this is under a different time parametrization. For each $t\ge0$, let
\begin{equation}
\tau_\varrho(t):=\inf\bigl\{s\ge0\colon\ell_s(\varrho)\ge t\bigr\}
\end{equation}
Then set
\begin{equation}
\label{E:1.31}
L_t(x):=\ell_{\tau_\varrho(t)}(x),\quad x\in D_N\cup\{\varrho\}
\end{equation}
Since $L_t(\varrho) = t$ a.s., this is the parametrization of the local time by the local time at~$\varrho$. 

As it turns out $E^\varrho(L_t(x))=t$ for each~$x\in D_N$, so \eqref{E:1.8} tells us that the cover time with happen on time scales such that $t\asymp (\log N)^2$. This motivates the parametrization \eqref{E:1.32} in our second main result to be discussed in these lectures:

\begin{theorem}[Abe-B.~\cite{AB}]
\label{thm-1.7}
Suppose $\{t_N\}_{N\ge1}$ is a positive sequence such that
\begin{equation}
\label{E:1.32}
\theta:=\lim_{N\to\infty}\frac{t_N}{2g(\log N)^2}>0
\end{equation}
Then setting
\begin{equation}
\label{E:1.33}
\wh K_N:=N^2\texte^{-\frac{t_N}{g\log N}}
\end{equation}
for $\theta\in(0,1)$ we have
\begin{equation}
\frac1{\wh K_N}\sum_{x\in D_N} 1_{\{L_{t_N}(x)=0\}}\delta_{x/N}\,\,\,\underset{N\to\infty}\lawarrow\,\,\,\,Z^D_{\sqrt\theta}
\end{equation}
where the random measures $\{Z_\lambda^D\colon\lambda\in(0,1)\}$ are as in Theorem~\ref{thm-1.5}. The limit (and, for large enough~$N$, the sequence of measures on the left) vanishes when~$\theta>1$.
\end{theorem}

Observe that from \twoeqref{E:1.32}{E:1.33} we get
\begin{equation}
\wh K_N = N^{2(1-\theta)+o(1)},\quad N\to\infty
\end{equation}
This decays to zero when~$\theta>1$ and so $\theta=1$ corresponds to the scaling of the cover time. (Note that we make no claims in this case as that requires going beyond the leading order asymptotic; see  Appendix of these notes for more discussion.)

The punchline of Theorem~\ref{thm-1.7} is that, in the parametrization by the local time at~$\varrho$, the set of avoided points at $\theta$-multiple of the cover time is asymptotically distributed \textit{exactly} as the $\sqrt\theta$-thick points of the DGFF. We note that, even though we are looking at the points where the local time vanishes, the time-parametrization matters. Indeed, a follow-up joint work with Y.~Abe and S.~Lee~\cite{ABL} shows that the limit measure is different in the natural time parametrization. See Theorem~\ref{thm-5.2} in Appendix for a precise statement.

\smallskip
As for the DGFF, we get information about the cardinality of the avoided set:

\begin{corollary}
\label{cor-1.8}
For the setting of Theorem~\ref{thm-1.7},
\begin{equation}
\frac1{\wh K_N}\#\bigl\{x\in D_N\colon L_{t_N}(x)=0\bigr\} \,\,\,\underset{N\to\infty}\lawarrow\,\,\,\,Z^D_{\sqrt\theta}(D)
\end{equation}
\end{corollary}

\smallskip
\noindent
Theorem~\ref{thm-1.7} will be extracted from Theorem~\ref{thm-1.5} and the connection of the local time and the DGFF discussed in Lecture~2. The actual proof will be given in Lecture~4.

\newpage
\section{Green function asymptotic and connection to DGFF}
\noindent
The purpose of the second lecture, which is also a kind of tutorial, is to develop two important technical ingredients that enter our later proofs. The first of these is the asymptotic of the Green function which is responsible for many of the underlying phenomena. The second ingredient concerns the connection of the local time to the DGFF which will later allow us to derive Theorem~\ref{thm-1.7} from Theorem~\ref{thm-1.5}.

\subsection{Asymptotic of the Green function}
Suppose that~$D$ is an admissible domain. For any~$x\in D$, let $\Pi^D$ denote the \textit{harmonic measure} from~$x$. This can be defined as the exit distribution from~$D$ of the standard Brownian motion~$B$ started at~$x$, i.e., for any Borel~$A\subseteq\R^2$,
\begin{equation}
\Pi^D(x,A):=P^x\bigl(B_{\tau_{D^\cc}}\in A\bigr)\quad\text{ where}\quad \tau_{D^\cc}:=\inf\bigl\{t\ge0\colon B_t\not\in D\bigr\}
\end{equation}
Let~$\{D_N\}_{N\ge1}$ be a sequence of admissible approximations of~$D$. Write $\lfloor xN\rfloor$ for the unique~$z\in\Z^d$ such that $x-z\in[0,1)^2$. We wish to prove:

\begin{theorem}
\label{thm-2.1}
Let~$\Vert\cdot\Vert$ denote the Euclidean norm on~$\R^2$. We then have:
\settowidth{\leftmargini}{(11)}
\begin{enumerate}
\item[(1)]
There exist a constant $c\in(0,\infty)$ such that for all $N\ge1$ and all $x,y\in D_N$
\begin{equation}
\label{E:2.2}
G^{D_N}(x,y)\le g\log\biggl(\frac{N}{1+\Vert x-y\Vert}\biggr)+c
\end{equation}
where~$g$ is as in \eqref{E:1.21}.
\item[(2)] For all $x\in D$,
\begin{equation}
\label{E:2.3}
G^{D_N}\bigl(\lfloor xN\rfloor,\lfloor xN\rfloor\bigr) = g\log N+c_0+g\log r^D(x)+o(1)
\end{equation}
where $c_0:=\frac14(2\gamma+\log 8)$ for~$\gamma$ denoting the Euler constant, 
\begin{equation}
\label{E:2.4}
r^D(x):=\exp\biggl\{\int_{\partial D}\Pi^D(x,\textd z)\log\Vert z-x\Vert\biggr\}
\end{equation}
and $o(1)\to0$ as~$N\to\infty$ locally uniformly in~$x\in D$.
\item[(3)] For all $x,y\in D$ with~$x\ne y$,
\begin{equation}
\label{E:2.5}
G^{D_N}\bigl(\lfloor xN\rfloor,\lfloor yN\rfloor\bigr)=-g\log\Vert x-y\Vert+g\int_{\partial D}\Pi^D(x,\textd z)\log\Vert z-y\Vert+o(1)
\end{equation}
where~$o(1)\to0$ as~$N\to\infty$ locally uniformly in~$(x,y)\in D\times D\smallsetminus\{(z,z)\colon z\in D\}$.
\end{enumerate}
\end{theorem}

Before we delve into the proof, let us make three remarks. First, the asymptotic statements \eqref{E:2.3} and \eqref{E:2.5} require that the points~$x$ and~$y$ stay away from~$\partial D$ and from each other. This is because the Green function vanishes near the boundary and has a logarithmic singularity on the ``diagonal.''

Second, the somewhat strange way of writing the~$x$ dependent term on the right of~\eqref{E:2.3} is motivated by the interpretation of~$r^D$. Indeed, for~$D$ simply connected this quantity coincides with the \textit{conformal radius} of~$D$ from~$x$, which is a measure of the distance of~$x$ to~$\partial D$ that is invariant under conformal maps. (The fact that $x\mapsto\log\Vert x-z\Vert$ is harmonic in~$x\ne z$ is used crucially in verifying this property.) In particular, $x\mapsto r^D(x)$ is a bounded continuous function on~$D$ tending to zero as~$x$ approaches~$\partial D$.

Third, the limit function on the right of \eqref{E:2.5}, namely,
\begin{equation}
\label{E:2.6}
\wh G^D(x,y):=-g\log\Vert x-y\Vert+g\int_{\partial D}\Pi^D(x,\textd z)\log\Vert z-y\Vert
\end{equation}
is the so called \textit{continuum Green function} in~$D$ with Dirichlet boundary condition. Since also this function is symmetric and positive semidefinite, it is a covariance, albeit only for a generalized Gaussian process called the \textit{Continuum Gaussian Free Field} (CGFF).

We remark that the CGFF is defined only by projections on suitable test functions (see, e.g., S.~Sheffield's review~\cite{Sheffield}) due to the fact that $\wh G^D(x,y)\to\infty$ as~$y\to x$ which makes pointwise value meaningless. This makes working with CGFF somewhat technically involved. Still, the singularity is only logarithmic and so thinking of the field as a random function usually gives a very good intuition.

\subsection{Proof of Theorem~\ref{thm-2.1}}
The proof of Theorem~\ref{thm-2.1} is based on a convenient representation of the Green function using the so called \textit{potential kernel}. In our normalization, this is a function~$\fraka\colon\Z^2\to\R$ defined by
\begin{equation}
\label{E:2.6a}
\fraka(x):=\frac14\int_{(-\pi,\pi)^2}\frac{\textd k}{(2\pi)^2}\,\frac{1-\cos(k\cdot x)}{\sin(k_1/2)^2+\sin(k_2/2)^2}
\end{equation}
where the integral converges because the numerator in the integrand vanishes quadratically (in~$k$) in the limit as~$k\to0$. With this we get:

\begin{lemma}
For all finite~$V\subseteq\Z^2$ and all~$x,y\in V$, 
\begin{equation}
\label{E:2.7}
G^V(x,y) = -\fraka(x-y)+\sum_{z\in\partial V}H^V(x,z)\fraka(z-y)
\end{equation}
where $H^V(x,z)$ is the probability that~$X$ started at~$x$ exists~$V$ at~$z\in\partial V$.
\end{lemma}

\begin{proofsect}{Proof}
Denote the discrete Laplacian acting on $f\colon\Z^2\to\R$ with compact support as
\begin{equation}
\Delta f(x):=\sum_{y\sim x}\bigl[f(y)-f(x)\bigr]
\end{equation}
where $y\sim x$ denotes that~$(x,y)$ is an edge in~$\Z^2$. Using the Markov property of~$X$ it is then  checked that, for each~$h\in V$, we have
\begin{equation}
\label{E:2.10i}
\left\{\begin{alignedat}{3}
\Delta G^V(\cdot,y)&=-\delta_y(\cdot)\quad &&\text{on }V
\\
G^V(\cdot,y)&=0\quad &&\text{on }\Z^2\smallsetminus V
\end{alignedat}\right.
\end{equation}
What makes the potential kernel particularly useful in this proof is that it solves a similar problem; namely,
\begin{equation}
\label{E:2.11i}
\left\{\begin{alignedat}{3}
\Delta \fraka(\cdot)&=\delta_0(\cdot)\quad &\text{on }\Z^2
\\
\fraka(0)&=0\quad &
\end{alignedat}\right.
\end{equation}
as is explicitly checked from \eqref{E:2.6a}.
Combining \twoeqref{E:2.10i}{E:2.11i} we conclude that
\begin{equation}
x\mapsto G^V(x,y)+\fraka(x-y)\text{ is discrete harmonic on~$V$}
\end{equation}
Relying on the fact that the discrete harmonic function is a martingale for the underlying random walk, we thus get
\begin{equation}
G^V(x,y) +\fraka(x-y)=\sum_{z\in\partial V}H^V(x,z)\bigl[G^V(z,y)+\fraka(z-y)\bigr]
\end{equation}
Noting that $G^V(z,y)=0$ for~$z\not\in V$, this reduces to \eqref{E:2.7}. 
\end{proofsect}

In order to make use of the formula \eqref{E:2.7} we need two lemmas whose proof we will leave to an exercise and/or literature study. The first of these concerns the asymptotic growth of the potential kernel, which is also where the constants~$g$ and~$c_0$ in Theorem~\ref{thm-2.1} enter the fray:

\begin{lemma}
\label{lemma-2.3}
For~$x\ne0$ we have $\fraka(x)\ge0$. Moreover,
\begin{equation}
\label{E:2.13}
\fraka(x) = g\log\Vert x\Vert+c_0+O\bigl(\Vert x\Vert^{-2}\bigr)
\end{equation}
\end{lemma}

This was apparently first proved by A.~St\"ohr~\cite{Stohr} in 1950. An article by G.~Kozma and E.~Schreiber~\cite{Kozma-Schreiber} from 2004 links the constant~$g$ and~$c_0$ to geometric properties of the underlying lattice, which allows then to verify the formula for other lattices as well. A very probabilistic approach to the theory of the potential kernel can be found in Section~4.4 of the monograph by G.~Lawler and V.~Limi\'c~\cite{Lawler-Limic}.

\begin{exercise}
Prove \eqref{E:2.13} by way of asymptotic analysis of the integral \eqref{E:2.6a}.
\end{exercise}

With Lemma~\ref{lemma-2.3} in hand, we are able to give:
\smallskip

\begin{proofsect}{Proof of (1) in Theorem~\ref{thm-2.1}}
Since $V\mapsto G^V(x,y)$ is non-decreasing with respect to the set inclusion, it suffices to prove this for~$x$ and~$y$ such that $d_\infty(y,D_N^\cc)\ge N$. Assuming~$x\ne y$ and plugging the asymptotic \eqref{E:2.13} shows
\begin{equation}
\begin{aligned}
G^{D_N}(x,y) = -\Bigl[\,g&\log\Vert x-y\Vert+c_0+O\bigl(\Vert x-y\Vert^{-2}\bigr)\Bigr]
\\&+\sum_{z\in\partial D_N}H^{D_N}(x,z)\Bigl[g\log\Vert z-y\Vert+c_0+O\bigl(\Vert z-y\Vert^{-2}\bigr)\Bigr]
\end{aligned}
\end{equation}
Using that $H^{D_N}(x,\cdot)$ is a probability mass function, the constant~$c_0$ cancels in both terms while the assumption that~$d_\infty(y,D_N^\cc)\ge N$ means that the last term in the square bracket is~$O(N^{-2})$. Bounding $\Vert z-y\Vert$ by a constant times~$N$, we then get the desired bound. In the case when~$x=y$ we use that~$\fraka(0)=0$ and bound the second term as above.
\end{proofsect}

For the remaining two parts of Theorem~\ref{thm-2.1} we also need:

\begin{lemma}
\label{lemma-2.5}
For any~$x\in D$,
\begin{equation}
\sum_{z\in\partial D_N} H^{D_N}\bigl(\lfloor xN\rfloor,z\bigr)\delta_{z/N}\,\,\overset{\text{\rm vaguely}}{\underset{N\to\infty}\longrightarrow}\,\, \Pi^D(x,\cdot)
\end{equation}
\end{lemma}

The proof of this is somewhat technical due to the fact that we want to make this work for rather general~$D$. The argument proceeds by coupling the random walk to Brownian motion so that their exit distributions remain close to each other. Details can be found in Appendix of a joint paper of the author with O.~Louidor~\cite{BL2}.

\smallskip
We are now ready for:

\begin{proofsect}{Proof of (2-3) in Theorem~\ref{thm-2.1}}
Starting with~(2), since~$\fraka(0)=0$ we only need to use the asymptotic on the second term on the right of \eqref{E:2.7}. This gives
\begin{equation}
\begin{aligned}
G^{D_N}\bigl(\lfloor xN\rfloor,&\lfloor xN\rfloor\bigr)=\sum_{z\in\partial D_N} H^{D_N}\bigl(\lfloor xN\rfloor,z\bigr)\fraka\bigl(z-\lfloor xN\rfloor\bigr)
\\
&=g\log N+c_0+g\sum_{z\in\partial D_N} H^{D_N}\bigl(\lfloor xN\rfloor,z\bigr)\log\biggl(\frac{\Vert z-\lfloor xN\rfloor\Vert}N\biggr)+O(N^{-2})
\end{aligned}
\end{equation}
where the error term uses the fact that, for~$N$ large enough, $\Vert z-z-\lfloor xN\rfloor\Vert\ge\delta N$ for some~$x$-dependent~$\delta>0$ uniformly in~$z\in\partial D_N$.
For similar reason the vague convergence in Lemma~\ref{lemma-2.5} applies to the function $z\mapsto\log\Vert z-x\Vert$ and, by way of an elementary approximation to get rid of the integer rounding, makes the sum to converge to~$\log r^D(x)$, locally uniformly in~$x\in D$.

For part~(3) we assume~$x\ne y$ and again use the asymptotic to get
\begin{equation}
\begin{aligned}
G^{D_N}\bigl(\lfloor xN\rfloor,\lfloor yN\rfloor\bigr)
=-&g\log\biggl(\frac{\Vert \lfloor xN\rfloor-\lfloor yN\rfloor\Vert}N\biggr)
\\
&+
g\sum_{z\in\partial D_N} H^{D_N}\bigl(\lfloor xN\rfloor,z\bigr)\log\biggl(\frac{\Vert z-\lfloor yN\rfloor\Vert}N\biggr)+O(N^{-2})
\end{aligned}
\end{equation}
where~$c_0$ again dropped out using that $z\mapsto H^{D_N}(x,z)$ is a probability mass function. Passing to the limit using Lemma~\ref{lemma-2.5} then yields the claim.
\end{proofsect}

\subsection{Connection between the local time and the DGFF}
The second topic of our interest in this lecture is a connection between the local time and the DGFF. We will treat this in the general case of a Markov chain on $V\cup\{\varrho\}$  where~$\varrho$ is the distinguished vertex that was used to define the Green function~$G^V$. 

Recall that~$L_t$ is the local time parametrized by the local time at~$\varrho$ which in particular means $L_t(\varrho)=t$ a.s. As our first result, we state the following limit theorem:

\begin{theorem}[DGFF limit]
\label{thm-2.6}
For~$L_t$ sampled under~$P^\varrho$,
\begin{equation}
\label{E:2.19}
\frac{L_t-t}{\sqrt{2t}}\,\,\,\underset{t\to\infty}\lawarrow\,\,\, h^V
\end{equation}
and, in particular,
\begin{equation}
\label{E:2.19a}
\sqrt{L_t}-\sqrt t\,\,\,\,\underset{t\to\infty}\lawarrow\,\,\, \frac1{\sqrt2} h^V
\end{equation}
where~$h^V$ is the DGFF on~$V$.
\end{theorem}

This result tells us that, at large times, a properly shifted and scaled local time profile is close to a sample of the DGFF. The rewrite \eqref{E:2.19a} explains why it is sometimes better with the square-root of~$L_t$ as no normalization is required. However, the connection runs far deeper and, in fact, applies at any fixed time~$t$:

\begin{theorem}[Second Ray-Knight Theorem]
\label{thm-2.7}
For each~$t\ge0$, there exists a coupling of~$L_t$ and two copies~$h$ and~$\tilde h$ of DGFF on~$V$ such that
\begin{equation}
\text{\rm $L_t$ and~$h$ are independent}
\end{equation}
and
\begin{equation}
\label{E:2.21}
\forall x\in V\cup\{\varrho\}\colon\quad L_t(x)+\frac12 h_x^2\,=\, \frac12\bigl(\tilde h_x+\sqrt{2t})^2\quad\text{\rm a.s.}
\end{equation}
\end{theorem}

\smallskip\noindent
This was proved as equality in distribution by N.~Eisenbaum, H.~Kaspi, M.B.~Marcus, J.~Rosen and Z.~Shi~\cite{EKMRS} in 2000 with the coupling part added by A.~Zhai~\cite{Zhai} in 2018.

We remark that Theorem~\ref{thm-2.7} belongs to a larger collection of results that link local time of stochastic processes to random fields. That such connection exists was first conceived of by K.~Symanzik~\cite{Symanzik}, and later developed by mathematical physicics (D.~Brydges, J.~Fr\"ohlich and T.~Spencer~\cite{BFS}) and probabilists (E.B.~Dynkin~\cite{Dynkin-book-chapter}). While Theorem~\ref{thm-2.7} is sometimes referred to as ``Dynkin isomorphism,'' this is a misattribution as the ``isomorphism'' in~\cite[Theorem~1]{Dynkin-book-chapter} works under a different setting than~\eqref{E:2.21}.

To explain the reliance on the specific time parameterization, 
observe that, under~$P^\varrho$, the local time $L_t$ on~$V$ is the sum of a random number of independent excursions that start by an exponential waiting time at~$\varrho$, then exit into~$V$ and, after running around~$V$ for a while, terminate by hitting~$\varrho$  again. Denoting, for each~$x\in V\cup\{\varrho\}$, the first return time of the chain to~$x$ by
\begin{equation}
\hat H_x:=\inf\bigl\{t\ge0\colon X_t=x\wedge\exists s\in[0,t)\colon X_s\ne x\bigr\}
\end{equation}
we thus get:

\begin{lemma}
\label{lemma-2.8}
Given $t>0$, let $\{\ell_j\}_{j\ge1}$ be i.i.d.\ copies of $\ell_{\hat H_\varrho}$ sampled under~$P^\varrho$ and let $N_t$ denote an independent Poisson random variable with parameter $\pi(\varrho)t$. Then
\begin{equation}
\label{E:2.23}
L_t \text{\rm\ under }P^\varrho\,\,\laweq \,\,\sum_{j=1}^{N_t}\ell_j, \quad\text{\rm on }V
\end{equation}
\end{lemma}

\begin{proofsect}{Proof}
The independence of the excursions is a consequence of the strong Markov property of~$X$. That the number of excursions has Poisson law is the standard fact that the number of i.i.d. Exponentials with parameter one  needed to accumulate the total value at least~$u$ is Poisson with parameter~$u$. 
\end{proofsect}

Note that \eqref{E:2.23} fails at~$\varrho$ a.s.\ due to the fact that~$L_t(\varrho)$ is non-random, which is the reason why we exclude~$\varrho$ from many statements below. With Lemma~\ref{lemma-2.8} in hand one can already prove the formula \eqref{E:2.19}. While we will  prove both theorems along the same lines, we leave the alternative argument to:

\begin{exercise}
Using the notation in Lemma~\ref{lemma-2.8}, prove that
\begin{equation}
\forall x\in V\cup\{\varrho\}\colon\,\,\,E^\varrho\bigl(\ell_1(x)\bigr)=1
\end{equation}
and
\begin{equation}
\forall x,y\in V\cup\{\varrho\}\colon\,\,\,\Cov_{P^\varrho}\bigl(\ell_1(x),\ell_1(y)\bigr)=G^V(x,y)
\end{equation}
Then use the multivariate (random-index) Central Limit Theorem and the decomposition in Lemma~\ref{lemma-2.8} to prove Theorem~\ref{thm-2.6}.
\end{exercise}

\subsection{Kac moment formula}
The proof of the above results is actually somewhat algebraic in nature. In order to present the details, introduce the standard inner product
\begin{equation}
\langle f,g\rangle:=\sum_{x\in V\cup\{\varrho\}}f(x)g(x), 
\end{equation}
and let $M_f$ be the operator of point-wise multiplication by~$f$ acting as
\begin{equation}
M_fg(x):=f(x)g(x), \quad x\in V\cup\{\varrho\}
\end{equation}
The driving force of all subsequent derivations is then:

\begin{lemma}[Kac moment formula] 
For each~$f\colon V\cup\{\varrho\}\to\R$ with~$f(\varrho)=0$,
\begin{equation}
\label{E:1.7}
E^\varrho\bigl(\langle \ell_1,f\rangle^n\bigr)=n!\frac1{\pi(\varrho)}\bigl\langle\, f,(G^VM_f)^{n-1}\,1\bigr\rangle,\quad n\ge1
\end{equation}
where~$(G^VM_f)g(x) = \sum_{y\in V\cup\{\varrho\}}G^V(x,y)f(y)g(y)$.
\end{lemma}

\begin{proofsect}{Proof}
The Markov property and elementary symmetrization tells us
\begin{equation}
\label{E:2.30}
E^\varrho\bigl(\langle \ell_1,f\rangle^n\bigr)=n!\int_{0\le t_1<\dots<t_n<H_\varrho}\textd t_1\dots\textd t_n\,\,\frac{f(X_{t_1})}{\pi(X_{t_1})}\dots \frac{f(X_{t_n})}{\pi(X_{t_n})}
\end{equation}
Abbreviating 
\begin{equation}
\cmss P^t(x,y):=\cmss P^x(X_t=y,\,\hat H_\varrho>t)
\end{equation}
and changing variables to $s_k:=t_k-t_{k-1}$ (where~$t_0:=0$), the Markov property of~$X$ allows us to rewrite the integral in \eqref{E:2.30} as
\begin{equation}
\label{E:1.12}
\sum_{x_1,\dots,x_n\in V\cup\{\varrho\}}\biggl(\prod_{i=1}^n\frac{f(x_i)}{\pi(x_i)}\biggr)\int_{s_1,\dots,s_n\ge0}\textd s_1\dots\textd s_n\,\,\cmss P^{s_1}(\varrho,x_1)\dots \cmss P^{s_n}(x_{n-1},x_n) 
\end{equation}
Next we observe that
\begin{equation}
\label{E:1.13a}
\int_0^\infty\textd s\,\,\cmss P^s(x,y)=E^x\Bigl(\,\int_0^{H_\varrho}\textd s\,1_{\{X_s=y\}}\Bigr)=\pi(y)G^V(x,y),\quad x,y\ne\varrho
\end{equation}
and, using the strong Markov property at the first hitting time of~$y$, 
\begin{equation}
\label{E:1.14a}
\begin{aligned}
\int_0^\infty\textd s\,\,\cmss P^s(\varrho,y)=P^\varrho(H_y<\hat H_\varrho)\pi(y)G^V(y,y)
,\quad x\ne \varrho
\end{aligned}
\end{equation}
To bring \eqref{E:1.14a} to a better form, use reversibility to get
\begin{equation}
\label{E:2.32}
\pi(\varrho)P^\varrho(H_y<\hat H_\varrho)=\pi(y)P^y(H_\varrho<\hat H_y)
\end{equation}
and then note that, by a decomposition of the form \eqref{E:1.13a} and the fact that Exponential(1)-random variable has mean one, $\pi(y)G^V(y,y)$ equals one plus the expected time to first succeed in independent trials with success probability $P^y(H_\varrho<\hat H_y)$. This implies
\begin{equation}
\label{E:2.33}
\pi(y)G^V(y,y)=\frac1{P^y(H_\varrho<\hat H_y)}
\end{equation}
which combining with \eqref{E:2.32} gives
\begin{equation}
\label{E:2.34}
\int_0^\infty\textd s\,\,\cmss P^s(\varrho,y)=\frac{\pi(y)}{\pi(\varrho)},\quad y\ne\varrho
\end{equation}
Note that this a different structure than \eqref{E:1.13a}.

For~$f\colon V\cup\{\varrho\}\to\R$ with~$f(\varrho)=0$ we now restrict the sums in \eqref{E:1.12} to~$x_i\in V$ and note that \eqref{E:1.13a} and \eqref{E:2.34} give
\begin{equation}
E^\varrho\bigl(\langle \ell_1,f\rangle^n\bigr)=\frac{n!}{\pi(\varrho)}\sum_{x_1,\dots,x_n\ne\varrho} f(x_1)G(x_1,x_2)\dots G(x_{n-1},x_n)f(x_n).
\end{equation}
The sum on the right is identified with $\langle f,(GM_f)^{n-1}1\rangle$.
\end{proofsect}

A formula of the kind \eqref{E:1.7} was first derived by M.~Kac~\cite{Kac} with a follow up by D.A.~Darling and M.~Kac~\cite{Darling-Kac}. A reader interested in more background and further results should consult P.J.~Fitzimmons and J.~Pitman~\cite{FP}.

\smallskip
The Kac moment formula gives us an explicit handle of the law of~$L_t$:

\begin{corollary}
For any $f\colon V\cup\{\varrho\}\to\R$ with $f(\varrho)=0$ and~$\max_{x\in V}|f(x)|$ small enough so that $\Vert G^VM_f\Vert<1$,
\begin{equation}
\label{E:1.16i}
E^\varrho\bigl(\texte^{\langle \ell_1,f\rangle}\bigr)=1+\frac1{\pi(\varrho)}\bigl\langle\, f,(1-G^VM_f)^{-1}1\bigr\rangle.
\end{equation}
In particular, for each~$t\ge0$,
\begin{equation}
\label{E:1.17i}
E^\varrho
\bigl(\texte^{\langle L_t,f\rangle}\bigr)=\texte^{\,t\langle\, f,(1-G^VM_f)^{-1}\,1\rangle}
\end{equation}
\end{corollary}

\begin{proofsect}{Proof}
Assume that~$f$ is so small that $\Vert GM_f\Vert<1$. The identity \eqref{E:1.16i} then follows by summing \eqref{E:1.7} on $n\ge1$. With the help from \eqref{E:2.23} we then get
\begin{equation}
E^\varrho
\bigl(\texte^{\langle L_t,f\rangle}\bigr)=\exp\Bigl\{t\pi(\varrho)\bigl[E\bigl(\texte^{\langle \ell_1,f\rangle}\bigr)-1\bigr]\Bigr\}
\end{equation}
and so \eqref{E:1.17i} follows from \eqref{E:1.16i}.
\end{proofsect}

This now permits us to conclude:

\begin{proofsect}{Proof of Theorem~\ref{thm-2.6}}
Assuming~$f$ small enough, rewrite \eqref{E:1.17i} as
\begin{equation}
\label{E:1.20i}
E^\varrho
\bigl(\texte^{\langle (L_t-t),f\rangle}\bigr)=\texte^{\,t\langle\, f,(1-G^VM_f)^{-1}\,G^V\,f\rangle}
\end{equation}
Now rescale~$f$ by~$\sqrt{2t}$ and notice that, as~$t\to\infty$, the right-hand side tends to $\texte^{\frac12\langle f,\,G^Vf\rangle}$, which is the Laplace transform of~$\langle f,h^V\rangle$. The Curtiss theorem then gives the claim.
\end{proofsect}

\subsection{Proof of the Second Ray-Knight Theorem}
Since~$V$ is fixed throughout, we will ease the notation by writing~$G$ for~$G^V$ throughout this subsection. In order to prepare for the proof of Theorem~\ref{thm-2.7}, we recall:

\begin{lemma}[Gaussian integration by parts]
Let $X=(X_1,\dots,X_n)$ be a multivariate Gaussian with mean zero and covariance matrix~$C$. Then for any $g\in C^1(\R^n)$ with subgaussian growth of~$\nabla g$ and any linear~$f\colon\R^n\to\R$,
\begin{equation}
\text{\rm Cov}\bigl(\,f(X),g(X)\bigr)=\sum_{i,j=1,\dots,n}C(i,j)E\biggl(\,\frac{\partial f}{\partial x_i}(X)\frac{\partial g}{\partial x_j}(X)\biggr)
\end{equation}
\end{lemma}

\begin{proofsect}{Proof}
For~$X_1,\dots,X_n$ i.i.d. $\NN(0,1)$, this is checked readily from $x\texte^{-\frac12x^2}=-\frac{\textd }{\textd x}\texte^{-\frac12x^2}$ and integration by parts. The general case is handled by writing $X=AZ$ where~$Z$ is a vector of i.i.d.~$\NN(0,1)$ and~$A$ is a matrix such that $\Cov(X) = AA^{\text{\rm T}}$.
\end{proofsect}

Using this we first note:

\begin{lemma}
For all $f,g\colon V\cup\{\varrho\}\to\R$ with~$f$ small enough and each~$s\in\R$,
\begin{equation}
\E\Bigl(\langle h+s,g\rangle\,\texte^{\frac12\langle (h+s)^2,f\rangle}\Bigr)
=s\bigl\langle 1,(1-M_f G)^{-1}g\bigr\rangle\, \E\Bigl(\texte^{\frac12\langle (h+s)^2,f\rangle}\Bigr)
\end{equation}
\end{lemma}

\begin{proofsect}{Proof}
Gaussian integration by parts shows
\begin{equation}
\begin{aligned}
\E\Bigl(\langle h+s,g\rangle\,&\texte^{\frac12\langle (h+s)^2,f\rangle}\Bigr)
=s\langle 1,g\rangle\,\E\Bigl(\texte^{\frac12\langle (h+s)^2,f\rangle}\Bigr)
+\E\Bigl(\langle h,g\rangle\,\texte^{\frac12\langle (h+s)^2,f\rangle}\Bigr)
\\
&=s\langle 1,g\rangle \E\Bigl(\texte^{\frac12\langle (h+s)^2,f\rangle}\Bigr)
+\E\Bigl(\bigl\langle M_f(h+s),Gg\bigr\rangle\,\texte^{\frac12\langle (h+s)^2,f\rangle}\Bigr)
\end{aligned}
\end{equation}
Putting the last term on the right together with the term on the left we get
\begin{equation}
\E\Bigl(\bigl\langle h+s,g-M_f Gg\bigr\rangle\,\texte^{\frac12\langle (h+s)^2,f\rangle}\Bigr)
=s\langle 1,g\rangle\, \E\Bigl(\texte^{\frac12\langle (h+s)^2,f\rangle}\Bigr)
\end{equation}
The claim follows by relabeling~$g$ for $(1-M_f G)^{-1}g$.
\end{proofsect}

Hence we get:

\begin{corollary}
For any $f\colon V\cup\{\varrho\}\to\R$ sufficiently small and any $s\ge t$,
\begin{equation}
\label{E:2.44}
\E\Bigl(\texte^{\frac12\langle (h+s)^2,f\rangle}\Bigr)
=\texte^{\frac12(s^2-t^2)\langle 1,(1-M_f G)^{-1}f\rangle}\,
\E\Bigl(\texte^{\frac12\langle (h+t)^2,f\rangle}\Bigr)
\end{equation}
\end{corollary}

\begin{proofsect}{Proof}
Using the previous lemma, we get
\begin{equation}
\begin{aligned}
\frac\textd{\textd r}
\E\Bigl(\texte^{\frac12\langle (h+r)^2,f\rangle}\Bigr)
&=\E\Bigl(\langle h+r,f\rangle\,\texte^{\frac12\langle (h+r)^2,f\rangle}\Bigr)
\\
&=r\bigl\langle 1,(1-M_f G)^{-1}f\bigr\rangle \,\E\Bigl(\texte^{\frac12\langle (h+r)^2,f\rangle}\Bigr)
\end{aligned}
\end{equation}
The differential equation is readily solved to get the result.
\end{proofsect}

\begin{proofsect}{Proof of Theorem~\ref{thm-2.7}}
Noting that
\begin{equation}
\texte^{\frac12(s^2-t^2)\langle 1,(1-M_f G)^{-1}f\rangle}
\end{equation}
is the exponential moment of $\langle L_r,f\rangle$ for $r := \frac12(s^2-t^2)$, we rewrite \eqref{E:2.44} as
\begin{equation}
\label{E:2.47}
\E\Bigl(\texte^{\frac12\langle (h+s)^2,f\rangle}\Bigr)
=E^\varrho\bigl(\texte^{\langle L_r,f\rangle}\bigr)\,
\E\Bigl(\texte^{\frac12\langle (h+t)^2,f\rangle}\Bigr)
\end{equation}
As this holds for all~$f$ small, solving for~$s$ as a function of~$t$ and~$r$ and applying the fact that the Laplace transform determines the underlying law shows
\begin{equation}
\label{E:2.48}
L_r\independent h\quad\Rightarrow\quad
L_r+\frac12(h+t)^2\,\,\laweq\,\,\frac12\bigl(h+\sqrt{t^2+2r}\,\bigr)^2
\end{equation}
Setting~$t:=0$ then gives \eqref{E:2.21} as equality in distribution. In order to construct the coupling, given independent~$L_t$ and~$h$, sample~$\tilde h$ from
\begin{equation}
\label{E:2.49}
\BbbP\Bigl(\, h^V\in\cdot\,\Big|\, \tfrac12\bigl(h^V+\sqrt{2t})^2=\phi\Bigr)\Big|_{\phi:=L_t+\frac12h^2}
\end{equation}
where the conditioning is well defined by the fact that the probability density of~$h^V$ is a continuous function.
The identity \eqref{E:2.21} then holds a.s.
\end{proofsect}

We finish with the following remark: Note that relabeling $t$ for~$\sqrt{2t}$ in \eqref{E:2.48} gives \eqref{E:2.21} in the form
\begin{equation}
\label{E:2.53}
L_r\independent h\quad\Rightarrow\quad L_r+\frac12\bigl(h^V+\sqrt{2t}\bigr)^2\,\,\laweq\,\,\frac12\Bigl(h^V+\sqrt{2(r+t)}\,\Bigr)^2
\end{equation}
which can alternatively be  derived by iterating \eqref{E:2.21} while using the independence of increments of~$t\mapsto L_t$. However, since the construction of the signs of $\tilde h+\sqrt{2t}$, which is what sampling from the conditional measure \eqref{E:2.49} is really about, is non-constructive, a question remains whether an almost-sure coupling can be constructed simultaneously for all times. We thus pose:

\begin{question}
Is there a coupling of the local time~$\{L_t\colon t\ge0\}$ (sampled under~$P^\varrho$) and an $\R^{V\cup\{\varrho\}}$-valued c\`adl\`ag process $\{h(t)\colon t\ge0\}$ such that
\begin{enumerate} 
\item[(1)] $\forall t\ge0\colon\,\, h(t)\,\laweq\,\,h^V$,
\item[(2)] $\forall t\ge0\colon\,\{h(s)\colon s\le t\}$ and $\{L_{t+u}-L_t\colon u\ge0\}$ are independent,
\item[(3)] for all $r,t\ge0$,
\begin{equation}
\label{E:2.55}
\forall r,t\ge0\colon\quad\frac12\bigl(h(r)+\sqrt{2r}\bigr)^2-L_r\,\,=\,\,\frac12\bigl(h(t)+\sqrt{2t}\,\bigr)^2 - L_t,\quad \text{a.s.}
\end{equation}
\end{enumerate}
hold true?
\end{question}

\noindent
To see that \eqref{E:2.55} is consistent with \eqref{E:2.53} note that, for $t\ge r$ we have $L_t-L_r\,\laweq L_{t-r}$ and so bringing $L_t$ to the left-hand side results in an identity that at least holds in distribution. 
The reason why we ask for such a coupling is two-fold. First, we find this to be an interesting possibility. Second, having the coupling would make some of the technical arguments in, e.g., \cite{ABL} much easier.

\newpage
\section{Thick points of the DGFF}
\noindent
The main goal of this lecture is to give the proof of Theorem~\ref{thm-1.5}. Not all of the details will be spelled out; the point is to convey the main ideas and explain the key technical steps. The reader is referred to the PIMS notes of the author~\cite{B-notes} for deeper treatment and, if even that is not sufficient, to the original joint paper with O.~Louidor~\cite{BL4}.

\subsection{Gibbs-Markov property of DGFF}
We start by an important fact about the Gaussian Free Field that we call the \textit{Gibbs-Markov property}. This is nothing but the ``domain-Markov property'' introduced for the continuum GFF in N.~Berestycki's lectures; the reason for attaching Gibbs' name to this concept is that, for DGFF, this property arises from the fact that the law of the DGFF is a Gibbs measure for a nearest-neighbor Hamiltonian. 
Here is the precise statement:

\begin{lemma}[Gibbs-Markov property]
\label{lemma-3.1}
Let~$U,V\subseteq\Z^2$ be non-empty finite sets with $U\subsetneq V$ and let~$h^V$ be the DGFF in~$V$. Recall that $H^U(x,z)$ is the probability that the simple random walk started at~$x$ exists~$U$ at~$z\in\partial U$ and set
\begin{equation}
\label{E:3.1}
\varphi^{V,U}_x:=\begin{cases}
\displaystyle\sum_{z\in\partial U}H^U(x,z)h^V_z,\qquad&\text{if }x\in U
\\*[5mm]
h_x^V,\qquad&\text{if }x\not\in U
\end{cases}
\end{equation}
Then 
\begin{equation}
\label{E:3.2}
h^V-\varphi^{V,U}\text{\rm\ and }\varphi^{V,U} \text{\rm\ are independent}
\end{equation}
with
\begin{equation}
\label{E:3.3}
h^V-\varphi^{V,U}\,\,\laweq \,\,\text{\rm DGFF in~$U$}
\end{equation}
Every sample path of~$\varphi^{V,U}$ is discrete harmonic on~$U$.
\end{lemma}

\begin{proofsect}{Proof}
If we set $H^U(x,z)=\delta_{xz}$ when~$x\not\in U$, we can write \eqref{E:3.1} concisely as
\begin{equation}
\label{E:3.4}
\varphi^{V,U}_x = \sum_{z\in V\smallsetminus U}H^U(x,z)h^V_z
\end{equation}
Noting that $H^U(x,\cdot)$ is a probability mass function, hence we get
\begin{equation}
\label{E:3.5}
\begin{aligned}
f(x,y):&=\Cov\bigl(h^V_x-\varphi^{V,U}_x,\varphi^{V,U}_y\bigr) 
\\
&= \sum_{z,z'\in V\smallsetminus U}\bigl[G^V(x,z')-G^V(z,z')\bigr]H^U(x,z)H^U(y,z')
\end{aligned}
\end{equation}
Now observe some facts about~$f$. First, $f(x,y)=0$ whenever~$x\in V\smallsetminus U$. Second, $x\mapsto f(x,y)$ is discrete harmonic on~$U$ and, for each $y\in V\smallsetminus U$, equals
\begin{equation}
\label{E:3.6}
G^V(x,y)-\sum_{z\in V\smallsetminus U}G^V(z,y)H^U(x,z)
\end{equation}
which is discrete harmonic in~$x\in U$ and vanishes at $x\in V\smallsetminus U$. The uniqueness of discrete harmonic extension forces~$f(\cdot,y)=0$ whenever~$y\in V\smallsetminus U$. From discrete harmonicity of~$y\mapsto f(x,y)$ we then get that~$f(x,y)=0$ for all~$x,y\in V$; i.e., 
\begin{equation}
\Cov\bigl(h^V_x-\varphi^{V,U}_x,\varphi^{V,U}_y\bigr)=0,\quad x,y\in V
\end{equation}
meaning that $h^V-\varphi^{V,U}$ and $\varphi^{V,U}$ are uncorrelated. As both fields are multivariate Gaussian, they are independent, proving \eqref{E:3.2}.

In order to prove \eqref{E:3.3} we use \eqref{E:3.4} to observe that
\begin{equation}
\Cov\bigl(h^V_x-\varphi^{V,U}_x,\,h^V_y-\varphi^{V,U}_y\bigr) = G^V(x,y)+\text{harmonic in~$x,y\in U$}
\end{equation}
The potential-kernel representation \eqref{E:2.7} then shows that, for each~$y\in U$,
\begin{equation}
\fraka(x-y)+\Cov\bigl(h^V_x-\varphi^{V,U}_x,\,h^V_y-\varphi^{V,U}_y\bigr)
\end{equation}
is harmonic in~$x\in U$ and equal to~$\fraka(z-y)$ at all~$z\in\partial U$. The argument in the proof of the identity \eqref{E:2.7} then gives
\begin{equation}
\Cov\bigl(h^V_x-\varphi^{V,U}_x,\,h^V_y-\varphi^{V,U}_y\bigr) = G^U(x,y)
\end{equation}
proving \eqref{E:3.3}. The discrete harmonicity of~$x\mapsto\varphi^{V,U}_x$ is a consequence of the same property of~$x\mapsto H^U(x,z)$ for~$z\not\in U$.
\end{proofsect}

\begin{remark}
Note that \eqref{E:3.1} means that~$\varphi^{V,U}$ is the harmonic extension of the values of~$h^V$ outside~$U$. Since $\varphi^{V,U}$ determines these values, \twoeqref{E:3.2}{E:3.3} in turn gives
\begin{equation}
\varphi^{V,U}_x = E\bigl(h^V_x\,\big|\,\sigma(h^V_z\colon z\in V\smallsetminus U)\bigr),\quad x\in V
\end{equation}
This is sometimes used as an alternative definition of~$\varphi^{V,U}$.
\end{remark}

In order to make referencing to the Gibbs-Markov property easier, we write it as
\begin{equation}
\label{E:3.11}
h^V\,\laweq\, h^U+\varphi^{V,U}\quad\text{where}\quad h^U\independent\varphi^{V,U}
\end{equation}
where~$h^V$ and~$h^U$ are DGFFs in~$V$ and~$U$ and~$\varphi^{V,U}$ has the law as specified above. 

One setting in which we will need the Gibbs-Markov property is when~$U:=V\smallsetminus\{x\}$. Then \eqref{E:3.11} reads
\begin{equation}
\label{E:3.12}
h^V\,\laweq\, \frakg_x(\cdot)h^V_x+h^U \quad\text{where}\quad h^U\independent h^V_x
\end{equation}
and~$\frakg_x\colon V\to[0,1]$ is a (deterministic) function that is discrete harmonic on~$V\smallsetminus\{x\}$ with $\frakg_x(x)=1$ and $\frakg_x=0$ on~$V^\cc$. As is easy to check, we have
\begin{equation}
\label{E:3.13}
\frakg_x(y) = \frac{G^V(x,y)}{G^V(x,x)},\quad y\in V
\end{equation}
This allows for control of~$\frakg_x$ using the estimates/asymptotics in Theorem~\ref{thm-2.1}.

\begin{figure}
\refstepcounter{obrazek}
\centerline{\includegraphics[width=1.9in]{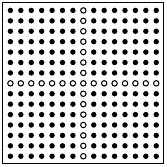}}
\begin{quote}
\label{fig6}
\fontsize{9}{5}\selectfont
{\bf Fig.~\theobrazek:\ }An illustration of the geometric setting for one typical use of the Gibbs-Markov property. Here $V$ is a box of~$(2N-1)^2$ vertices which is split into four $(N-1)\times(N-1)$ squares (whose union is~$U$) and a ``cross'' made of two lines of vertices separating these.
\normalsize
\medskip
\end{quote}
\end{figure}

\begin{figure}[t]
\refstepcounter{obrazek}
\centerline{\includegraphics[width=3.6in]{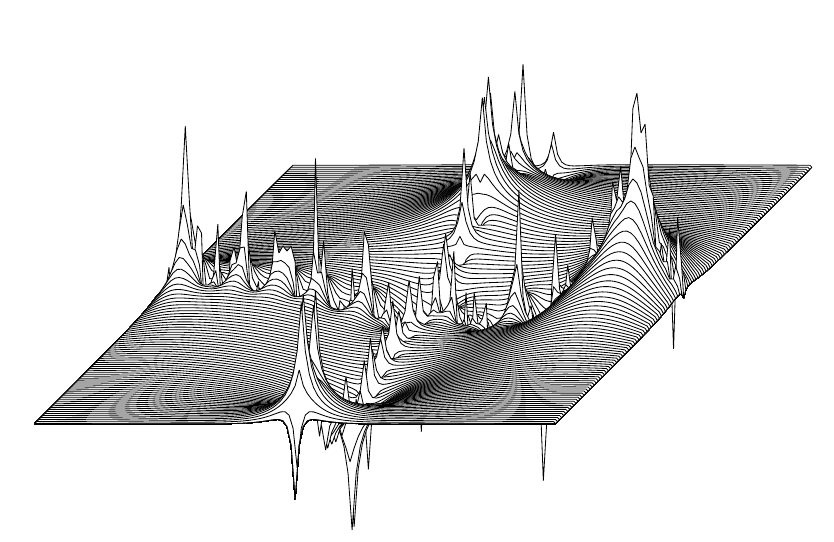}} 
\begin{quote}
\label{fig8}
\fontsize{9}{5}\selectfont
{\bf Fig.~\theobrazek:\ }A sample of $\varphi^{V,U}$ for the geometric setting in Fig.~\ref{fig6}. The field is discrete harmonic and thus ``smooth'' on~$U$ but is quite rough on the ``cross'' of vertices constituting~$V\smallsetminus U$, due to the fact that it coincides with~$h^V$ there.
\normalsize
\medskip
\end{quote}
\end{figure}

Another instance where the Gibbs-Markov property will be used is when~$V$ is the square $\{1,\dots,2N-1\}^2$ and~$U$ is the union of four translates of the square $\{1,\dots,N-1\}^2$ by vectors $(0,0)$, $(N,0)$, $(0,N)$ and~$(N,N)$. The set~$V\smallsetminus U$ is a ``cross'' made of two lines of vertices; see Fig.~\ref{fig6}. Note that, on the ``cross,'' $\varphi^{V,U}$ has the law of~$h^V$ and is thus quite rough there. However, thanks to discrete harmonicity, $\varphi^{V,U}$ is quite smooth once sufficiently far from the boundary. A sample of $\varphi^{V,U}$ is shown in Fig.~\ref{fig8}.

\smallskip

An important fact associated with the Gibbs-Markov property in domains $U\subseteq V$ that are scaled-up version of two continuum domains by~$N$, the field~$\varphi^{V,U}$ is well approximated by, and in fact converges to, a smooth process. Explicitly, we have:

\begin{lemma}
\label{lemma-3.2}
Let $\{D_N\}_{N\ge1}$ and~$\{\wt D_N\}_{N\ge1}$ be admissible approximation of two admissible domains~$\wt D\subseteq D\subseteq\R^2$. For each~$N\ge1$ there exists a coupling of $\varphi^{D_N,\wt D_N}$ and a Gaussian process $\{\Phi^{D,\wt D}(x)\colon x\in \wt D\}$ with law determined by
\begin{equation}
\forall x\in \wt D\colon\quad\E\Phi^{D,\wt D}(x)=0
\end{equation}
and
\begin{equation}
\forall x,y\in \wt D\colon\quad \E\bigl(\Phi^{D,\wt D}(x)\Phi^{D,\wt D}(y)\bigr)=\wh G^D(x,y)-\wh G^{\wt D}(x,y),
\end{equation}
where~$\wh G^D$ is the continuum Green function in~$D$ defined in \eqref{E:2.6},
such that for each $\delta>0$,
\begin{equation}
\sup_{\begin{subarray}{c}
x\in \wt D\\ d(x,\wt D^\cc)>\delta
\end{subarray}}
\bigl|\varphi^{D_N,\wt D_N}_{\lfloor xN\rfloor}-\Phi^{D,\wt D}(x)\bigr|\,\,\underset{N\to\infty}{\overset{P}\longrightarrow}\,\,0
\end{equation}
Moreover, a.e.\ sample path of~$\Phi^{D,\wt D}$ is harmonic on~$\wt D$.
\end{lemma}

Note that the singular parts of the continuum Green function cancel in the expression $\wh G^D(x,y)-\wh G^{\wt D}(x,y)$ and so the expression is meaningful even when~$x=y$. That this is a covariance kernel (with respect to the Lebesgue measure) follows from it being the limit of the covariances of $\varphi^{D_N,\wt D_N}$. This limit implies weak convergence in the sense of finite-dimensional distributions. To get convergence in supremum norm, one has to control the oscillation of the two processes which is done using concentration techniques for Gaussian processes; see, e.g.,~\cite[Chapter~6]{B-notes}.

\smallskip

\subsection{Subsequential limits}
We now move to the main goal of this section. Fix a sequence $\{a_N\}_{N\ge1}$ such that~$\lambda$ defined by the limit in \eqref{E:1.25} belongs to~$(0,1)$. We will only carry out the proof under the assumption that~$\lambda<1/\sqrt2$ because this does not require truncations in second-moment calculations we perform below.

Let~$\{D_N\}_{N\ge1}$ be admissible approximations of an admissible domain~$D\subseteq\R^2$. Abbrevite the measures of interest as
\begin{equation}
\label{E:3.18i}
\eta_N:=\frac1{K_N}\sum_{x\in D_N} \delta_{x/N}\otimes\delta_{h^{D_N}_x-a_N}
\end{equation}
and recall that $K_N=N^{2(1-\lambda^2)+o(1)}$. We also introduce the notation
\begin{equation}
\Gamma_N(b):=\bigl\{x\in D_N\colon h_x^{D_N}\ge a_N+b\bigr\}
\end{equation}
for the level set of~$h^{D_N}$ at ``height'' $a_N+b$. Using the shorthand, $\langle\mu,f\rangle:=\int f\textd\mu$, note that, given any $A\subseteq D$ and abbreviating $A_N:=\{x\in\Z^2\colon x/N\in A\}$, we have
\begin{equation}
\frac1{K_N}\bigl|\Gamma_N(b)\cap A_N\bigr|=\bigl\langle \eta_N,1_A\otimes1_{[b,\infty)}\bigr\rangle
\end{equation}
Our next goal is to compute two moments of the size of the level set $\Gamma_N(b)$. We start with a bound on the first moment:

\begin{lemma}
\label{lemma-3.3}
There exists~$c>0$ such that for all $N\ge1$ and all $b\in[-\frac12a_N,a_N]$,
\begin{equation}
\forall A\subseteq D_N\colon\quad \E\bigl|\Gamma_N(b)\cap A\bigr|\le c \frac{|A|}{N^2} \,\texte^{-\frac{a_N}{g\log N}b}K_N
\end{equation}
\end{lemma}

\begin{proofsect}{Proof}
Recall that $G(x,x):=G^{D_N}(x,x)\le g\log N+\tilde c$ uniformly in~$x\in D_N$. Using the standard Gaussian estimate
\begin{equation}
\label{E:Gauss-bd}
X=\NN(0,\sigma^2)\quad\Rightarrow\quad\forall t\ge0\colon\,\,P(X\ge t)\le \sigma t^{-1}\texte^{-\frac{t^2}{2\sigma^2}} 
\end{equation}
we get
\begin{equation}
\begin{aligned}
\BbbP\bigl(h_x^{D_N}\ge a_N+b\bigr)
&\le\frac{\sqrt{G(x,x)}}{a_N+b}\,\texte^{-\frac{(a_N+b)^2}{2G(x,x)}}
\\
&\le \frac{\sqrt{g\log N+\tilde c}}{a_N/2}\,\texte^{-\frac{(a_N+b)^2}{2[g\log N+\tilde c]}}
\le \frac c{\sqrt{\log N}}\,\texte^{-\frac{a_N^2}{2g\log N}}\,\texte^{-\frac{a_N}{g\log N} b}
\end{aligned}
\end{equation}
where we used that $(g\log N+\tilde c)^{-1} = (g\log N)^{-1}+O(\log N)^{-2})$. Invoking the definition of~$K_N$, the claim follows by summing over~$x\in A$.
\end{proofsect}

The purpose of the bound in Lemma~\ref{lemma-3.3} is to give us control over the contribution from the part of~$D_N$ close to the boundary where the asymptotic for the Green function from Theorem~\ref{thm-2.1} does not apply. This helps to get:

\begin{lemma}
\label{lemma-3.4}
Let~$A\subseteq D$ be open and set $A_N:=\{x\in\Z^2\colon x/N\in A\}$. Then for all~$b\in\R$,
\begin{equation}
\frac1{K_N}\E\bigl|\Gamma_N(b)\cap A_N\bigr|\,\,\underset{N\to\infty}\longrightarrow\,\,\hat c(\alpha\lambda)^{-1}\texte^{-\alpha\lambda b}\int_A r^D(x)^{2\lambda^2}\textd x
\end{equation}
where $\alpha:=2/\sqrt g$ and
\begin{equation}
\hat c:=\frac{\texte^{-2c_0\lambda^2/g}}{\sqrt{2\pi g}}
\end{equation}
for~$c_0$ the constant in \eqref{E:2.3} and~$r^D$ is as in \eqref{E:2.4}.
\end{lemma}

\begin{proofsect}{Proof}
Lemma~\ref{lemma-3.3} allows us to assume that~$A_N$ is at least~$\delta N$ away from~$D_N^\cc$ and focus on~$x\in A_N$ with~$h^{D_N}_x\le a_N+\log\log N$. The asymptotic in Theorem~\ref{thm-2.1}(2) gives
\begin{equation}
\begin{aligned}
\BbbP\bigl(a_N+&\,\log\log N\ge h_x^{D_N}\ge a_N+b\bigr)
\\
&=\bigl(1+o(1)\bigr)\int_b^{\log\log N}\frac1{\sqrt{2\pi g\log N}}\,\,\texte^{-\frac{a_N^2}{2[g\log N+c_0+g\log r^D(x/N)]}}\,\texte^{-s\frac{a_N}{g\log N}}\
\textd s
\end{aligned}
\end{equation}
Using that $a_N/g\log N = \alpha\lambda+o(1)$ we can write the integrand as
\begin{equation}
\frac{1+o(1)}{\sqrt{2\pi g}}
\texte^{-2c_0\lambda^2/g}\,r^D(x/N)^{2\lambda^2}\texte^{-\alpha\lambda s}\,\frac{K_N}{N^2}
\end{equation}
Integrating over~$s$ and turning the sum over~$x\in A_N$ normalized by~$N^2$ into an integral then yields the claim.
\end{proofsect}

We now move to the second moment of~$|\Gamma_N(b)|$. It is here where the restriction on~$\lambda$ simplifies calculations quite a bit.

\begin{lemma}
\label{lemma-3.5}
Suppose~$\lambda\in(0,1/\sqrt2)$. For all~$b\in\R$ ther exists~$c>0$ such that for all~$N\ge1$,
\begin{equation}
\E\Bigl(\bigl|\Gamma_N(b)\bigr|^2\Bigr)\le c K_N^2
\end{equation}
\end{lemma}

\begin{proofsect}{Proof}
Let us set~$b:=0$ for simplicity (or absorb the term into~$a_N$). Then 
\begin{equation}
\E\Bigl(\bigl|\Gamma_N(b)\bigr|^2\Bigr) = \sum_{x,y\in D_N}\BbbP\bigl(h^{D_N}_x\ge a_N,\,h^{D_N}_y\ge a_N\bigr)
\end{equation}
To bound the summand uniformly in~$x$ or~$y$ regardless how far these are from the boundary~$D_N$, we replace~$h^{D_N}$ by the field in the enlarged domain
\begin{equation}
\wt D_N:=\{x\in\Z^2\colon d_\infty(x,D_N)\le N\}
\end{equation}
noting that then
\begin{equation}
\label{E:3.6a}
\frac14\BbbP\bigl(h^{D_N}_x\ge a_N,\,h^{D_N}_y\ge a_N\bigr)\le \BbbP\bigl(h^{\wt D_N}_x\ge a_N,\,h^{\wt D_N}_y\ge a_N\bigr)
\end{equation}
holds by a routine application of the  Gibbs-Markov property. 

Next we split the sum according to whether~$d_\infty(x,y)\le\sqrt{K_N}$ or not. The first part we bound using Lemma~\ref{lemma-3.3} as
\begin{equation}
\sum_{\begin{subarray}{c}
x,y\in D_N\\ d_\infty(x,y)\le\sqrt{K_N}
\end{subarray}}
\!\!\!\!\BbbP(h^{\wt D_N}_x\ge a_N,\,h^{\wt D_N}_y\ge a_N)
\le(2\sqrt{K_N}+1)^2\sum_{x\in D_N}\BbbP(h^{\wt D_N}_x\ge a_N)\le c K_N^2
\end{equation}
In the second part we distinguish whether $h^{\wt D_N}_x$ exceeds~$2a_N$ or not. 
Using that
\begin{equation}
\BbbP(h_x^{\wt D_N}\ge 2a_N\bigr)\le\frac{c}{\sqrt{\log N}}\texte^{-2\frac{a_N^2}{g\log N}}
=c\Bigl(\frac{K_N}{N^2}\Bigr)^2\sqrt{\log N}\,\texte^{-\frac{a_N^2}{g\log N}}
\le c \Bigl(\frac{K_N}{N^2}\Bigr)^2
\end{equation}
once~$N$ is sufficiently large, the part where $h^{\wt D_N}_x\ge 2a_N$ contributes at most
\begin{equation}
\sum_{x,y\in D_N}\BbbP\bigl(h^{\wt D_N}_x> 2a_N,\,h^{\wt D_N}_y\ge a_N\bigr)
\le c \Bigl(\frac{|D_N|}{N^2}\Bigr)^2 K_N^2
\end{equation}
where the right-hand side is again at most a constant times~$K_N^2$.

We are thus left to bound the expression
\begin{equation}
\label{E:3.9a}
\sum_{\begin{subarray}{c}
x,y\in D_N\\ d_\infty(x,y)>\sqrt{K_N}
\end{subarray}}
\BbbP\bigl(2a_N\ge h^{\wt D_N}_x\ge a_N,\,h^{\wt D_N}_y\ge a_N\bigr)
\end{equation}
Here we note that the Gibbs-Markov property allows us to condition on $h^{\wt D_N}_x\ge a_N$ by way of the decomposition \eqref{E:3.12} that in the present setting reads
\begin{equation}
h_y^{\wt D_N} \laweq \frakg_x(y) h_x^{\wt D_N} +\hat h_y^{\wt D_N\smallsetminus\{x\}}\end{equation}
where~$\frakg_x\colon\Z^2\to[0,1]$ is the unique discrete-harmonic function in $D_N\smallsetminus\{x\}$ extending the boundary values $\frakg_x(x)=1$\, and~$\,\frakg_x=0$ on~$D_N^\cc$
and $\hat h_y^{\wt D_N\smallsetminus\{x\}}$ is the DGFF in~$D_N\smallsetminus\{x\}$ that is sampled independently of $h_x^{\wt D_N}$. Using this we can write
\begin{equation}
\label{E:3.12a}
\begin{aligned}
\BbbP\bigl(2a_N\ge&\, h^{\wt D_N}_x\ge a_N,\,h^{\wt D_N}_y\ge a_N\bigr) 
\\
&= \int_{0}^{a_N}\BbbP\bigl(h_x^{\wt D_N}-a_N\in\textd s\bigr)\BbbP\Bigl(\hat h_y^{\wt D_N\smallsetminus\{x\}}\ge a_N\bigl[1-\frakg_x(y)\bigr]-s\frakg_x(y)\Bigr)
\end{aligned}
\end{equation}
In order to bound the integrand, observe that $d_\infty(x,y)>\sqrt{K_N}=N^{1-\lambda^2+o(1)}$ along with the fact that~$x,y$ are ``deep'' inside~$\wt D_N$ enable Theorem~\ref{thm-2.1} to give
\begin{equation}
\label{E:3.13a}
\begin{aligned}
\frakg_x(y)=\frac{G^{\wt D_N}(x,y)}{G^{\wt D_N}(x,x)}&\le\frac{\log\frac N{\Vert x-y\Vert}+c}{\log N-c}
\\
&\le 1-(1-\lambda^2)+o(1)=\lambda^2+o(1)
\end{aligned}
\end{equation}
Since $\lambda<1/\sqrt2$, it follows that, given any ~$\epsilon\in(0,1-2\lambda^2)$,
\begin{equation}
\epsilon a_N\le a_N\bigl[1-\frakg_x(y)\bigr]-s\frakg_x(y)\le a_N
\end{equation}
holds for all for $s\in[0,a_N]$  once~$N$ is sufficiently large, uniformly in~$x$ and~$y$ contributing to the sum \eqref{E:3.9a}. The standard Gaussian estimate \eqref{E:Gauss-bd}
along with the bound
\begin{equation}
\bigl(a_N[1-\frakg_x(y)]-s\frakg_x(y)\bigr)^2\ge a_N^2-2(a_N+s)a_N\frakg_x(y)
\end{equation}
 then show
\begin{equation}
\begin{aligned}
\BbbP\Bigl(\hat h_y^{\wt D_N\smallsetminus\{x\}}\ge&\, a_N\bigl[1-\frakg_x(y)\bigr]-s\frakg_x(y)\Bigr)
\\
&\le\frac{\sqrt{G(y,y)}}{\epsilon a_N}\,\texte^{-\frac{(a_N[1-\frakg_x(y)]-s\frakg_x(y))^2}{2G(y,y)}}
\le c\frac{K_N}{N^2}\,\texte^{\frac{a_N^2}{g\log N}\frakg_x(y)+\frac{a_N}{g\log N}\frakg_x(y)s}
\end{aligned}
\end{equation}
where~$G(y,y)$ abbreviates $G^{\wt D_N\smallsetminus\{x\}}(y,y)$ and where we also used $G(y,y)\le g\log N+c$ uniformly in~$y\in\wt D_N\smallsetminus\{x\}$.

Now observe that the first inequality in \eqref{E:3.13a} gives
\begin{equation}
\texte^{\frac{a_N^2}{g\log N}\frakg_x(y)}\le c\biggl(\frac N{\Vert x-y\Vert}\biggr)^{4\lambda^2+o(1)}
\end{equation}
Using the explicit form of the probability density of $h^{\wt D_N}_x$ we also get
\begin{equation}
\BbbP\bigl(h_x^{\wt D_N}-a_N\in\textd s\bigr)
\le c\frac{K_N}{N^2}\texte^{-\frac{a_N}{g\log N}s}\textd s
\end{equation}
With the help of these we bring \eqref{E:3.12a} to the form
\begin{equation}
\begin{aligned}
\BbbP\bigl(2a_N\ge&\, h^{\wt D_N}_x\ge a_N,\,h^{\wt D_N}_y\ge a_N\bigr) 
\\
&\le c\biggl(\frac{K_N}{N^2}\biggr)^2
\biggl(\frac N{\Vert x-y\Vert}\biggr)^{4\lambda^2+o(1)}
\int_{0}^{a_N}
\texte^{-\frac{a_N}{g\log N}[1-\frakg_x(y)]s}\textd s
\end{aligned}
\end{equation}
In light of the uniform bound $\frakg_x(y)\le\lambda^2+o(1)$, the integral converges uniformly in all~$y$ of concern. As a consequence, we get
\begin{equation}
\begin{aligned}
\sum_{\begin{subarray}{c}
x,y\in D_N\\ d_\infty(x,y)>\sqrt{K_N}
\end{subarray}}
\BbbP\bigl(2a_N\ge h^{\wt D_N}_x\ge& a_N,\,h^{\wt D_N}_y\ge a_N\bigr)
\\
&\le c\biggl(\frac{K_N}{N^2}\biggr)^2\!\!\!
\sum_{\begin{subarray}{c}
x,y\in D_N\\ d_\infty(x,y)>\sqrt{K_N}
\end{subarray}}\!\!\biggl(\frac N{\Vert x-y\Vert}\biggr)^{4\lambda^2+o(1)}
\end{aligned}
\end{equation}
Using that~$4\lambda^2<2$, the sum is dominated by pairs~$x$ and~$y$ such that~$\Vert x-y\Vert$ is order~$N$ and so is order~$(N^2)^2$. (Alternatively, dominate the sum by an integral.) The expression is thus bounded by a constant times~$K_N^2$, as desired.
\end{proofsect}

\subsection{Factorization and uniqueness}
As a consequence of Lemmas~\ref{lemma-3.3}--\ref{lemma-3.5} we get:

\begin{corollary}
\label{cor-3.6}
Suppose~$\lambda<1/\sqrt2$. Then~$\{\eta_N\}_{N\ge1}$ defined in \eqref{E:3.18i} is a tight sequence of measures on~$\overline D\times(\R\cup\{+\infty\})$ and every subsequential weak limit~$\eta$ obeys
\begin{equation}
\label{E:3.43}
\E\,\eta\bigl(A\times[b,\infty)\bigr) = \hat c(\alpha\lambda)^{-1}\texte^{-\alpha\lambda b}\int_A r^D(x)^{2\lambda^2}\textd x
\end{equation}
for any open~$A\subseteq D$ and any~$b\in\R$.
\end{corollary}

A particular consequence of \eqref{E:3.43} is that $\eta$ is non-vanishing on each non-empty open set with positive probability. Similarly as Lemma~\ref{lemma-3.4} refines the bound from Lemma~\ref{lemma-3.3}, the second moment calculation from Lemma~\ref{lemma-3.5} can be refined to get:

\begin{lemma}
\label{lemma-3.7}
Suppose~$\lambda<1/\sqrt2$ and, given~$A\subseteq D$ open, set $A_N:=\{x\in\Z^2\colon x/N\in A\}$. Then for all~$b\in\R$,
\begin{equation}
\frac1{K_N^2}\E\biggl(\Bigl[\bigl|\Gamma_N(b)\cap A_N\bigr|-\texte^{-\alpha\lambda b}\bigl|\Gamma_N(0)\cap A_N\bigr|\Bigr|^2\biggr)\,\,\underset{N\to\infty}\longrightarrow\,\,0
\end{equation}
\end{lemma}

\begin{proofsect}{Proof (idea)}
We write the expectation as the sum over~$x,y\in A_N$ of
\begin{equation}
\begin{aligned}
\BbbP\bigl(h_x^{D_N}&\ge a_N+b, h_y^{D_N}\ge a_N+b\bigr)
-\texte^{-\alpha\lambda b}\BbbP\bigl(h_x^{D_N}\ge a_N+b, h_y^{D_N}\ge a_N\bigr)
\\&-\texte^{-\alpha\lambda b}\BbbP\bigl(h_x^{D_N}\ge a_N, h_y^{D_N}\ge a_N+b\bigr)
+\texte^{-2\alpha\lambda b}\BbbP\bigl(h_x^{D_N}\ge a_N, h_y^{D_N}\ge a_N\bigr)
\end{aligned}
\end{equation}
The proof of Lemma~\ref{lemma-3.5} (with $\sqrt{K_N}$ cut-off replaced by $\delta\sqrt{K_N}$) tells us that it suffices to control the pairs $d_\infty(x,y)\ge \delta N$. We then need to show that for any $b_1,b_2\in\{0,b\}$,
\begin{equation}
\begin{aligned}
\BbbP\bigl(h_x^{D_N}&\ge a_N+b, h_y^{D_N}\ge a_N+b\bigr)
\\&=\bigl(\texte^{-\alpha\lambda(b_1+b_2)}+o(1)\bigr)\BbbP\bigl(h_x^{D_N}\ge a_N, h_y^{D_N}\ge a_N\bigr)
\end{aligned}
\end{equation}
which is checked by similar calculations as those in the proof of Lemmas~\ref{lemma-3.4} and~\ref{lemma-3.5}.
\end{proofsect}

We now upgrade the conclusion of Corollary~\ref{cor-3.6} as:

\begin{corollary}[Factorization]
Suppose~$\lambda<1/\sqrt2$. Then every subsequential weak limit~$\eta$ of measures $\{\eta_N\}_{N\ge1}$ factors as
\begin{equation}
\label{E:3.47}
\eta(\textd x\textd h) = Z^D(\textd x)\otimes \texte^{-\alpha\lambda h}\textd h
\end{equation}
where $Z^D$ is a random Borel measure such that
\begin{equation}
\label{E:3.48}
\E Z^D(A) = \hat c\int_A r^D(x)^{2\lambda^2}\textd x
\end{equation}
holds for any open~$A\subseteq D$. In particular, $Z^D(A)=0$ a.s.\ for any~$A\subseteq\R^d$ with vanishing Lebesgue measure.
\end{corollary}

\begin{proofsect}{Proof}
Lemma~\ref{lemma-3.7} shows that, for all $b\in\R$ and all~$A\subseteq D$ open,
\begin{equation}
\label{E:3.49}
\eta\bigl(A\times[b,\infty)\bigr) = \texte^{-\alpha\lambda b}\eta\bigl(A\times[0,\infty)\bigr)\text{\rm\ \ a.s.}
\end{equation}
Density arguments from measure theory permit us to choose the same implicit null set for all~$b$ and~$A$. Setting
\begin{equation}
Z^D(A):=\alpha\lambda\,\eta\bigl(A\times[0,\infty)\bigr)
\end{equation}
 the right hand side of \eqref{E:3.49} then coincides with the integral of the measure on the right of \eqref{E:3.47} over $A\times[b,\infty)$. The identity \eqref{E:3.48} then follows from \eqref{E:3.43}.
\end{proofsect}

With all the subsequential limit measures taking the desired product form, the following questions remain: Is the subsequential weak limit unique in law? And, if so, is there a way to characterize~$Z^D$? In order to answer these, we need to elucidate how the Gibbs-Markov property manifests itself for the limit object:

\begin{lemma}
\label{lemma-3.10}
Let~$D,\wt D$ be admissible domains with $\wt D\subseteq D$ yet with $D\smallsetminus \wt D$ of vanishing Lebesgue measure. Then for $Z^D$ and~$Z^{\wt D}$ constructed along the same subsequence,
\begin{equation}
\label{E:3.51}
Z^D(\textd x)\,\laweq\,\texte^{\alpha\lambda\Phi^{D,\wt D}(x)}Z^{\wt D}(\textd x),\quad \Phi^{D,\wt D}\independent Z^{\wt D}
\end{equation}
where on the right~$\Phi^{D,\wt D}$ is the Gaussian process from Lemma~\ref{lemma-3.2}. 
\end{lemma}

\begin{proofsect}{Proof}
Let $\eta^D_N$, resp, $\eta^{\wt D}_N$ denote the finite~$N$ processes in $D_N$, resp., $\wt D_N$. The Gibbs-Markov property in Lemma~\ref{lemma-3.1} gives a coupling of $\eta^D_N$, $\eta^{\wt D}_N$ and~$\varphi^{D_N,\wt D_N}$ such that
\begin{equation}
\bigl\langle \eta_N^D, f(\cdot,\cdot)\bigr\rangle = \bigl\langle \eta_N^{\wt D}, f(\cdot,\cdot+\varphi_{\lfloor \cdot N\rfloor}^{D_N,\wt D_N})\bigr\rangle\text{\rm\ \ a.s.},\quad \varphi^{D_N,\wt D_N}\independent \eta^{\wt D}_N
\end{equation}
whenever~$f$ is continuous and compactly supported in~$\wt D\times\R$. Passing to a joint distributional limit along the same subsequence yields a coupling of the limiting processes $\eta^D$ and~$\eta^{\wt D}$ and the field~$\Phi^{D,\wt D}$ such that $ \Phi^{D,\wt D}$ and $\eta^{\wt D}$ are independent and 
\begin{equation}
\langle \eta^D,f\rangle = \bigl\langle\eta^{\wt D},f(\cdot,\cdot+\Phi^{D,\wt D})\bigr\rangle\text{\rm\ \ a.s.}
\end{equation}
By separability arguments, the a.s.-equality can be made to work for all~$f$ simultaneously.
Plugging in \eqref{E:3.47}, a change of coordinates shows
\begin{equation}
Z^D(\textd x) = \texte^{\alpha\lambda\Phi^{D,\wt D}(x)}Z^{\wt D}(\textd x)\text{\rm\ \ a.s.}
\end{equation}
as measures on~$\wt D$. Since neither sides charges $D\smallsetminus \wt D$, the a.s.-equality applies even as measures on~$D$.
\end{proofsect}

We are finally ready to give:

\begin{proofsect}{Proof of Theorem~\ref{thm-1.5} for $\lambda<1/\sqrt2$, modulo a technical step}
We will argue that the law of $Z^D$ is uniquely determined by the expectation \eqref{E:3.48} and the Gibbs-Markov property \eqref{E:3.51}. For this we first extract a joint subsequential weak limit~$\eta^D$ for all~$D$ ranging through finite unions of dyadic squares in~$\R^2$ and define~$Z^D$ as above. 

Denote $S:=(0,1)\times(0,1)$. Use $\{S^n_i\colon i=1,\dots,4^n\}$ to label the squares of the form $(k2^{-n},\ell2^{-n})+(0,2^{-n})\times(0,2^{-n})$ for $k,\ell=0,\dots,2^{n}-1$ and denote
\begin{equation}
S^n:=\bigcup_{i=1}^{4^n}S^n_i
\end{equation}
Note that~$S\smallsetminus S^n$ has vanishing Lebesgue measure.
Lemma~\ref{lemma-3.10} then gives
\begin{equation}
\label{E:3.55}
Z^S(\textd x)\,\,\laweq\,\,\texte^{\alpha\lambda\Phi^{S,S^n}(x)}Z^{S^n}(\textd x),\quad \Phi^{S,S^n}\independent Z^{S^n}
\end{equation}
Since the DGFF is independent over the connected components of the underlying domain, the measures $\{Z^{S^n_i}\colon i=1,\dots, 4^n\}$ are independent and
\begin{equation}
\label{E:3.56}
\texte^{\alpha\lambda\Phi^{S,S^n}(x)}Z^{S^n}(\textd x) = \sum_{i=1}^{4^n}1_{S^n_i}(x) \texte^{\alpha\lambda\Phi^{S,S^n}(x)} Z^{S^n_i}(\textd x)
\end{equation}
Introduce the measure in which each $Z^{S^n_i}$ is replaced by its expectation,
\begin{equation}
Y^S_n(\textd x):=\hat c\sum_{i=1}^{4^n}1_{S^n_i}(x) \texte^{\alpha\lambda\Phi^{S,S^n}(x)} r^{S^n_i}(x)^{2\lambda^2}\textd x
\end{equation}
Let~$f\colon S\to[0,\infty)$ be continuous with compact support. Rewriting $\langle Z^S,f\rangle$ via \twoeqref{E:3.55}{E:3.56} and taking the conditional expectation over~$\{Z^{S^n_i}\colon i=1,\dots,4^n\}$ given~$\Phi^{S,S^n}$ with the help of \eqref{E:3.48} and the conditional Jensen inequality gives
\begin{equation}
\E\bigl(\texte^{-\langle Z^S,f\rangle}\bigr)\ge \E\bigl(\texte^{-\langle Y^S_n,f\rangle}\bigr)
\end{equation}
Our goal for much of the rest of the proof is to show the reverse inequality, at least in the limit as~$n\to\infty$. First we observe:

\begin{exercise}[Reverse Jensen inequality]
Prove that if $X_1,\dots,X_n$ are independent non-negative random variables, then
\begin{equation}
\label{E:3.59}
E\biggl(\exp\Bigl\{-\sum_{i=1}^n X_i\Bigr\}\biggr)\le \exp\biggl\{-\texte^{-\epsilon}\sum_{i=1}^n E(X_i|X_i\le\epsilon)\biggr\}
\end{equation}
holds for all $\epsilon>0$.
\end{exercise}

Given $\delta\in(0,1/2)$ we apply this to
\begin{equation}
X_i:=\bigl\langle Z^{S^n_i},\, 1_{S^{n,\delta}_i}f\bigr\rangle
\end{equation}
where $S^{n,\delta}_i$ is obtained by the same translate as $S^n_i$ but of the box $(\delta 2^{-n},(1-\delta)2^{-n})$. Since
\begin{equation}
\langle Z^S,f\rangle\ge\sum_{i=1}^{4^n}X_i
\end{equation}
the inequality \eqref{E:3.59} gives
\begin{equation}
\E\bigl(\texte^{-\langle Z^S,f\rangle}\bigr)\le 
\E\biggl(\exp\Bigl\{-\texte^{-\epsilon}\sum_{i=1}^{4^n} E(X_i|X_i\le\epsilon)\Bigr\}\biggr)
\end{equation}
where $E(X_i|X_i\le\epsilon)$ is still a conditional expectation given~$\Phi^{S,S_n}$.

In order to show that the approximations involving~$\epsilon$ and~$\delta$ are negligible, we then have to first show that the errors incurred by conditioning are negligible,
\begin{equation}
\forall\epsilon>0\colon\quad\sum_{i=1}^{4^n} E(X_i|X_i>\epsilon)\,\,\underset{n\to\infty}{\overset{P}\longrightarrow}\,\,0
\end{equation}
This is done by using Markov's inequality to write $E(X_i|X_i>\epsilon)\le\epsilon^{-1}E(X_i^2)$ and estimating the result using calculations similar to those in the proof of Lemma~\ref{lemma-3.5}. (Here assuming $\lambda<1/\sqrt2$ again makes life easier.)

Second, denoting $S^n_\delta:=\bigcup_{i=1}^{4^n}S^{n,\delta}_i$, we have to show that the restriction of the integral represented by $\langle Y^S_n,f\rangle$ from $S$ to $S^n_\delta$ is negligible,
\begin{equation}
\forall\epsilon>0\colon\quad\lim_{\delta\downarrow0}\limsup_{n\to\infty} \BbbP\bigl(Y^S_n(S\smallsetminus S^n_\delta)>\epsilon\bigr)=0
\end{equation}
This is easier as the first moment of $Y^S_n(S\smallsetminus S^n_\delta)$ is controlled by the Lebesgue measure of $S\smallsetminus S^n_\delta$. We refer to \cite[pages~217-218]{B-notes} for full details.

Since convergence of Laplace transforms of integrals of random measures with respect to compactly-supported continuous functions implies weak convergence of the measures, we have shown
\begin{equation}
\label{E:3.65}
Y^S_n\,\,\underset{n\to\infty}\lawarrow\,\, Z^S
\end{equation}
This holds regardless of the subsequence leading to the definition of~$Z^S$ and so we can now conclude the claim for~$D$ being a unit square. An analogous proof (or a scaling relation) extends this to all dyadic squares. Invoking Lemma~\ref{lemma-3.10} one more time along partitions of~$D$ into disjoint open dyadic squares then proves convergence for all admissible~$D$. (With the limit identified uniquely, we can also write~$Z^D_\lambda$.)
\end{proofsect}

We remark that \eqref{E:3.65} gives~$Z^S_\lambda$ a form of Gaussian multiplicative chaos associated with the continuum GFF which is central to N.~Berestycki's lectures. In that subject area, the measure~$Z^D_\lambda$ is referred as the \textit{Liouville quantum gravity}. The upcoming book by N.~Berestycki and E.~Powell~\cite{Berestycki-Powell} treats this subject area from a number of angles. See Fig.~\ref{fig9} for a sample of~$Z^D_\lambda$-measure.

\refstepcounter{obrazek}
\centerline{\includegraphics[width=3.6in]{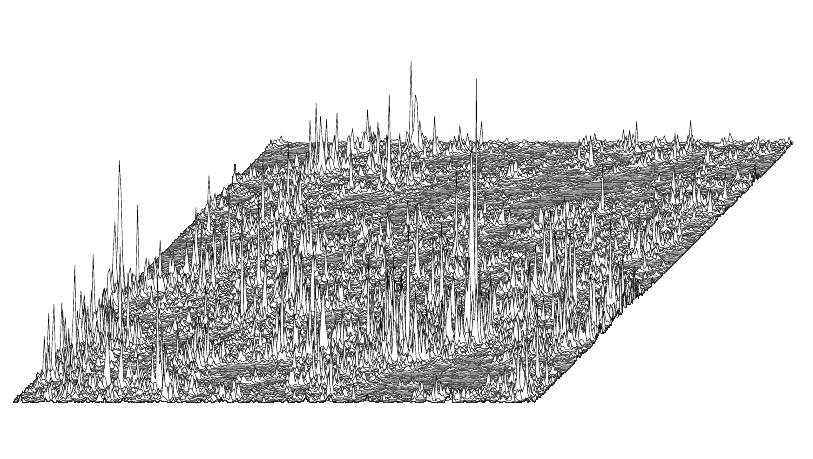}} 
\begin{quote}
\label{fig9}
\fontsize{9}{5}\selectfont
{\bf Fig.~\theobrazek:\ }A sample of $Z^D_\lambda$ for~$D:=(0,1)\times(0,1)$ and~$\lambda:=0.3$. The measure is known to be supported on a set of Hausdorff dimension $2(1-\lambda^2)$ in accord with $K_N=N^{2(1-\lambda^2)+o(1)}$.
\normalsize
\medskip
\end{quote}

Another point to note is that
\begin{equation}
\Var\bigl(\Phi^{D,\wt D}(x)\bigr) = g\log\biggl(\frac{r^D(x)}{r^{\wt D}(x)}\biggr)
\end{equation}
along with the fact (implied by $\Cov(\Phi^{D,\wt D})=\wh G^D-\wh G^{\wt D}$) that the increment fields
\begin{equation}
\{\Phi^{S,S^{n+1}}-\Phi^{S,S^n}\colon n\ge0\}
\end{equation}
can be realized as independent on the same probability space shows that $\langle Y^S_n,f\rangle$ is, for any continuous compactly-supported $f\ge0$, a positive martingale. The convergence \eqref{E:3.65} thus takes place a.s. Since $\wh G^D-\wh G^{\wt D}$ are all non-negative, a criterion of J.P.~Kahane~\cite{Kahane} then characterizes the law of~$Z^S_\lambda$ uniquely. The Gibbs-Markov property then extends this to~$Z^D_\lambda$ for all admissible~$D$. We refer to~\cite[Section~5.2]{B-notes} for details.

\newpage
\section{Points avoided by random walk}
\noindent
The fourth and final lecture of this minicourse is devoted to the proof of Theorem~\ref{thm-1.7} describing the scaling limit of the set of points avoided by the simple random walk on~$D_N\cup\{\varrho\}$. For reasons explained then, we will work with the time parametrization by the local time at the boundary vertex~$\varrho$.

\subsection{Setting the scales and main idea}
Let us start by explaining the formula for the normalizing sequence~$\{\wh K_N\}_{N\ge1}$ from \eqref{E:1.33}. Working for a moment in the general setting of a Markov chain on~$V\cup\{\varrho\}$, recall the definition \eqref{E:1.31} of the local time~$L_t$ parametrized by the local time at~$\varrho$. We then have:

\begin{lemma}
\label{lemma-4.1}
For all~$x\in V$ and all~$t\ge0$,
\begin{equation}
P^\varrho\bigl(L_t(x)=0\bigr) = \texte^{-\frac{t}{G^V(x,x)}}
\end{equation}
\end{lemma}

\begin{proofsect}{Proof}
Let~$x\in V$. Observe that, by the (a.s.-unique) time~$s$ when~$\ell_s(\varrho)=t$, the chain started at~$\varrho$ has accumulated a Poisson($\pi(\varrho)t$) number of excursions into~$V$. The probability that such an excursion visits~$x$ is $P^\varrho(H_x<\wh H_\varrho)$, which by a Poisson-thinning argument means that the total number of excursions that visit~$x$ by the time when the time at~$\varrho$ equals~$\pi(\varrho)t$ has the law of
\begin{equation}
\text{Poisson}\Bigl(\pi(\varrho)P^\varrho(H_x<\wh H_\varrho)t\Bigr)
\end{equation}
If $L_t(x)=0$, then no such excursion has visited~$x$ and so 
\begin{equation}
P^\varrho\bigl(L_t(x)=0\bigr) = \texte^{-\pi(\varrho)P^\varrho(H_x<\wh H_\varrho)t}
\end{equation}
We now observe that, by \twoeqref{E:2.32}{E:2.33}, the exponent equals $t/G^V(x,x)$.
\end{proofsect}

For~$V:=D_N$ we have $G^{D_N}(x,x)\le g\log N+c$ for~$x$ sufficiently far away from the boundary and so, for~$t=O((\log N)^2)$,
\begin{equation}
\label{E:4.4}
P^\varrho\bigl(\exists x\in D_N\colon L_t(x)=0\bigr)\le
E^\varrho\biggl(\,\sum_{x\in D_N}1_{\{L_t(x)=0\}}\biggr)\le c |D_N|\texte^{-\frac{t}{g\log N}}
\end{equation}
Hereby we conclude:

\begin{corollary}
For~$t_N\sim 2g\theta(\log N)^2$ with~$\theta>1$, we have
\begin{equation}
P^\varrho\bigl(\forall x\in D_N\colon L_t(x)>0\bigr)\,\,\underset{N\to\infty}\longrightarrow\,\,1
\end{equation}
\end{corollary}

\noindent
For~$\theta<1$ the calculation \eqref{E:4.4} shows that the expected number of avoided points grows as a power of~$N$. As we will show by proving Theorem~\ref{thm-1.7}, the actual number of avoided points runs on the same scale. In particular, we get that $\theta=1$ markes the leading-order time scale of the cover time.

As discussed in Lecture~2, our key tool will be the Second Ray-Knight Theorem (Theorem~\ref{thm-2.7}) that states that there exists a coupling of~$L_t$ and two copies ~$h$ and~$\tilde h$ of the DGFF on~$V$ such that
\begin{equation}
\label{E:4.5}
L_t\independent h\quad\wedge\quad L_t+\frac12 h^2\,=\, \frac12\bigl(\tilde h+\sqrt{2t}\bigr)^2
\end{equation}
We will use this roughly as follows: Suppose~$x$ is such that~$L_t(x)=0$. Then, if~$h_x$ happens to be order unity, also the field on the right is of order unity which means that
\begin{equation}
\tilde h_x = -\sqrt{2t}+O(1)
\end{equation}
For the choice~$t=2g\theta(\log N)^2$ we have $\sqrt{2t} = 2\sqrt g\sqrt\theta\log N$, which means that~$x$ is a $\sqrt\theta$-thick point of~$\tilde h$!

Of course, assuming that~$h_x$ is order unity needs to be justified because $h_x$ is (at typical~$x$) normal with variance~$\log N$ and so it is order unity only with probability proportional to $(\log N)^{-1/2}$. Since~$h$ is independent and the points where~$L_t$ vanishes are somewhat scattered, we interpret the above argument though Poisson thinning: Diluting the points where~$L_t(x)=0$ (more or less) independently with probability of order $(\log N)^{-1/2}$ gives us, roughly, the $\sqrt\theta$-points of~$\tilde h$.

\subsection{Light points}
In order to implement the above strategy quantitatively, a moment's thought reveals that tracking only the points where the local time vanishes is not sufficient. Instead, we will need to track the set where the local time is as well of order unity; we refer to these as the \textit{light points}. We introduce the notation for the corresponding measure
\begin{equation}
\vartheta_N:=\frac1{\wh K_N}\sum_{x\in D_N}\delta_{x/N}\otimes\delta_{L_t(x)}
\end{equation}
Our first goal is to show that the family of measures $\{\vartheta_N\}_{N\ge1}$ is tight. For this we need to upgrade Lemma~\ref{lemma-4.1} to the form:

\begin{lemma}
\label{lemma-4.3}
For all~$x\in V$, all $t\ge0$ and all~$b\ge0$,
\begin{equation}
P^\varrho\bigl(L_t(x)\le b\bigr) \le \texte^{-\frac{t}{G^V(x,x)}\exp\{-\frac{b}{G^V(x,x)}\}}
\le\texte^{-\frac{t}{G^V(x,x)}+b\frac{t}{G^V(x,x)^2}}
\end{equation}
\end{lemma}

\begin{proofsect}{Proof}
The argument from the proof of Lemma~\ref{lemma-4.1} gives us the representation
\begin{equation}
\label{E:4.10}
\pi(x) L_t(x)\,\laweq\,\sum_{k=1}^{N_t}\sum_{j=1}^{Z_k}T_{k,j}
\end{equation}
where
\begin{itemize}
\item $\{T_{k,j}\}_{k,j\ge1}$ are i.i.d.\ Exponential with parameter 1
\item $\{Z_k\}_{k\ge1}$ are i.i.d.\ Geometric with parameter $P^x(H_\varrho<\wh H_x)$
\item $N_t$ is Poisson with parameter $\pi(\varrho)P^\varrho(H_x<\wh H_\varrho)t$.
\end{itemize}
with all the random variables independent. Indeed, all we need to realize that, when an excursion from~$\varrho$ hits~$x$, the number of visits (i.e., the arrival plus returns) to~$x$ on this excursion will be Geometric with parameter $P^x(H_\varrho<\wh H_x)$. Each visits adds an independent Exponential(1)-time to the time spent at~$x$.

Abbreviate $q:=P^x(H_\varrho<\wh H_x)$. A thinning argument for exponential random variables gives
\begin{equation}
\sum_{j=1}^{Z_k}T_{k,j}\,\,\laweq\,\,\text{Exponential}(q)
\end{equation}
Observe also that reversibility \eqref{E:2.32} gives
\begin{equation}
\pi(\varrho)P^\varrho(H_x<\wh H_\varrho) = \pi(x)P^x(H_\varrho<\wh H_x)=\pi(x)q
\end{equation}
and so $N_t=\text{Poisson}(t\pi(x)q)$. Hence we get
\begin{equation}
\label{E:4.13}
\begin{aligned}
P^\varrho\bigl(L_t(x)\le b\bigr) &\le P\Biggl(\forall k=1,\dots, N_t\colon \sum_{j=1}^{Z_k}T_{k,j}\le \pi(x) b\Biggr)
\\
&=\sum_{n=0}^\infty \frac{[t\pi(x)q]^n}{n!}\bigl[1-\texte^{-\pi(x)b q}\bigr]^n\texte^{-t\pi(x)q}=\texte^{-t\pi(x)q\exp\{b\pi(x)q\}}
\end{aligned}
\end{equation}
To get the first bound in the claim, we observe that $\pi(x)q = G^V(x,x)^{-1}$. The second bound follows from the inequality $\texte^{-s}\ge 1-s$.
\end{proofsect}

Hereby we conclude:

\begin{corollary}
\label{cor-4.4}
Suppose~$t_N\sim 2g\theta(\log N)^2$ for~$\theta\in(0,1)$. For all~$b\ge0$ there exists~$c\in(0,\infty)$ such that for all~$A\subseteq \R^2$ open and with~$A_N:=\{x\in D_N\colon x/N\in A\}$,
\begin{equation}
\label{E:4.15}
E^\varrho\,\vartheta\bigl(A\times[0,b]\bigr)\le c\frac{|A_N|}{N^2}
\end{equation}
holds. In particular, $\{\vartheta_N\}_{N\ge1}$ are tight as measures on~$\overline D\times\R_+$.
\end{corollary}

\begin{proofsect}{Proof}
Let~$b\ge0$. Writing $G(x,x)$ instead of~$G^{D_N}(x,x)$, Lemma~\ref{lemma-4.3} gives
\begin{equation}
E^\varrho\,\vartheta\bigl(A\times[0,b]\bigr) \le\frac1{\wh K_N} \sum_{x\in A_N}
\min\Bigl\{\texte^{-\frac{t_N}{G(x,x)}+b\frac{t_N}{G(x,x)^2}},\,\texte^{-\frac{t_N}{G(x,x)}\exp\{-\frac{b}{G(x,x)}\}}\Bigr\}
\end{equation}
For~$x$ such that~$G(x,x)\ge \texte^{-b/4} g\log N$, the uniform bound~$G(x,x)\le g\log N+c$ implies
\begin{equation}
\frac{t_N}{G(x,x)}-b\frac{t_N}{G(x,x)^2}\ge \frac{t_N}{g\log N+c}-b\texte^{b/2}\frac{t_N}{(g\log N)^2}
\end{equation}
which is at least $\frac{t_N}{g\log N}-c'$ for some constant~$c'$ depending only on~$b$. If~$x$ in turn satisfies $G(x,x)\le \texte^{-b/4} g\log N$ then the fact that $G(x,x)\ge1/4$ implies
\begin{equation}
\frac{t_N}{G(x,x)}\exp\Bigl\{-\frac{b}{G(x,x)}\Bigr\}\ge \texte^{b/4}\frac{t_N}{g\log N}\,\texte^{-b/4}=\frac{t_N}{g\log N}
\end{equation}
Hence the minimum in \eqref{E:4.15} is at most a constant times $\texte^{-\frac{t_N}{g\log N}}=\wh K_N/N^2$. 
\end{proofsect}

\newcommand{\ext}{{\text{\rm ext}}}
\newcommand{\Leb}{{\text{\rm Leb}}}

\subsection{Extended process}
The tightness of $\{\vartheta_N\}_{N\ge1}$ permits us to extract subsequential weak limits. Our goal is to characterize these limits via the coupling \eqref{E:4.5} of~$L_t$, $h$ and~$\tilde h$ which in turn requires that, along with small values of~$L_t$ we also track small values of~$h$. For this we introduce the extended process
\begin{equation}
\vartheta^\ext_N:=\frac{\sqrt{\log N}}{\wh K_N}\sum_{x\in D_N}\delta_{x/N}\otimes\delta_{L_t(x)}\otimes\delta_{h_x}
\end{equation}
where the additional $\sqrt{\log N}$ in the normalization reflects on the fact that, forcing a point with a small value of~$L_t(x)$ to have a small value of~$h_x$ costs $O((\log N)^{-1/2})$ in probability. A key technical lemma to prove is:

\begin{lemma}
\label{lemma-4.5}
Suppose $\{N_k\}_{k\ge1}$ is a strictly increasing sequence such that $\vartheta_{N_k}\lawarrow\,\vartheta$. Then
\begin{equation}
\label{E:4.19}
\vartheta^\ext_{N_k}\,\,\underset{k\to\infty}\lawarrow\,\,\vartheta\otimes\Leb
\end{equation}
where~$\Leb$ is the Lebesgue measure on~$\R$.
\end{lemma}

\begin{proofsect}{Proof (modulo a technical step)}
We first note that the convergence takes place under expectation. Indeed, writing~$\E$ for the expectation with respect to~$h$ only (i.e., conditional on~$L_t$), for any $f=f(x,\ell,h)$ non-negative  and continuous with compact support,
\begin{equation}
\label{E:4.20}
\begin{aligned}
\E\bigl(\langle\vartheta^\ext_N, f\rangle\bigr) &= \frac{\sqrt{\log N}}{\wh K_N}\sum_{x\in D_N}\E f\bigl(x/N,L_t(x),h_x\bigr)
\\
&=\frac{\sqrt{\log N}}{\wh K_N}\sum_{x\in D_N}\int\frac1{\sqrt{2\pi}}\frac1{\sqrt{G(x,x)}}\,\texte^{-\frac{h^2}{2G(x,x)}} f\bigl(x/N,L_t(x),h\bigr)\textd h
\end{aligned}
\end{equation}
where we again write~$G(x,x)$ instead of~$G^{D_N}(x,x)$.
The restriction on support of~$f$ means that the integral over~$h$ is over a bounded interval and the sum over~$x$ at least~$\delta N$ away from the boundary of~$D_N$. This means that
\begin{equation}
\frac1{\sqrt{G(x,x)}}\,\texte^{-\frac{h^2}{2G(x,x)}} = \frac1{\sqrt{g\log N}}+O\Bigl(\frac1{\log N}\Bigr) 
\end{equation}
uniformly in~$h$ and~$x$ of interest. Since $2\pi g=1$,  we thus get
\begin{equation}
\label{E:4.22}
\E\bigl(\langle\vartheta^\ext_N, f\rangle\bigr)=\biggl(1+O\Bigl(\frac1{\sqrt{\log N}}\Bigr)\biggr)\langle\vartheta_N\otimes\Leb, f\rangle
\end{equation}
which, in light of tightness of~$\{\vartheta_N\}_{N\ge1}$, equals $\langle\vartheta_N\otimes\Leb f\rangle+o(1)$.

Next we observe that the conditional Jensen inequality upgrades the above to the form
\begin{equation}
\label{E:4.23}
E^\varrho\otimes \E\bigl(\texte^{-\langle\vartheta_N^\ext,f\rangle}\bigr)
\ge\texte^{o(1)} E^\varrho\bigl(\texte^{-\langle\vartheta_N\otimes\Leb,f\rangle}\bigr)
\end{equation}
Since convergence of Laplace transforms implies convergence in law, we thus need to show that the opposite inequality holds as~$N\to\infty$. This is done roughly as follows: Given~$s>0$, consider the expectation
\begin{equation}
E^\varrho\otimes \E\bigl(\langle\vartheta_N^\ext,f\rangle\texte^{-s\langle\vartheta_N^\ext,f\rangle}\bigr)
\end{equation}
This is related to the above by
\begin{equation}
\label{E:4.25}
E^\varrho\otimes \E\bigl(\texte^{-\langle\vartheta_N^\ext,f\rangle}\bigr)
=1-\int_0^s E^\varrho\otimes \E\bigl(\langle\vartheta_N^\ext,f\rangle\,\texte^{-s\langle\vartheta_N^\ext,f\rangle}\bigr)\textd s
\end{equation}
The additive form of $\langle\vartheta_N^\ext,f\rangle$ now implies
\begin{equation}
\label{E:4.26}
\begin{aligned}
E^\varrho\otimes \E\bigl(&\langle\vartheta_N^\ext,f\rangle\,\texte^{-s\langle\vartheta_N^\ext,f\rangle}\bigr)
\\
&=\frac{\sqrt{\log N}}{\wh K_N}\sum_{x\in D_N}
E^\varrho\otimes \E\Bigl(f\bigl(x/N,L_t(x),h_x\bigr)\texte^{-s\langle\vartheta_N^\ext,f\rangle}\Bigr)
\end{aligned}
\end{equation}
Since $f\ge0$ and the argument of~$f$ under expectation depends only on~$h_x$, the conditional Jensen inequality gives
\begin{equation}
\begin{aligned}
E^\varrho\otimes \E\Bigl(f\bigl(x/N,L_t(x),h_x\bigr)&\,\texte^{-s\langle\vartheta_N^\ext,f\rangle}\Bigr)
\\&\ge 
E^\varrho\otimes \E\Bigl(f\bigl(x/N,L_t(x),h_x\bigr)\,\texte^{-s\E(\langle\vartheta_N^\ext,f\rangle\,|\,h_x)}\Bigr)
\end{aligned}
\end{equation}
The point is to show that the conditioning on~$h_x$ can be ignored and the conditional expectation can be replaced by~$\langle\vartheta_N\otimes\Leb,f\rangle+o(1)$. This is done through a truncation argument for which we refer the reader to \cite[Lemma~7.1]{AB}.

Once the conditioning is taken care of, we again apply the calculation \twoeqref{E:4.20}{E:4.22} to the remaining occurrence of~$h_x$ in the expression under the sum in \eqref{E:4.26}. Hence we get
\begin{equation}
E^\varrho\otimes \E\bigl(\langle\vartheta_N^\ext,f\rangle\,\texte^{-s\langle\vartheta_N^\ext,f\rangle}\bigr)
\ge o(1)+\texte^{o(1)}
E^\varrho\bigl(\langle\vartheta_N\otimes\Leb,f\rangle\,\texte^{-s\langle\vartheta_N\otimes\Leb,f\rangle}\bigr)
\end{equation}
where both~$o(1)$ tend to zero as~$N\to\infty$ uniformly in~$s\in[0,1]$.
Plugging this in \eqref{E:4.25} then gives
\begin{equation}
E^\varrho\otimes \E\bigl(\texte^{-\langle\vartheta_N^\ext,f\rangle}\bigr)
\le o(1)+\texte^{o(1)} E^\varrho\bigl(\texte^{-\langle\vartheta_N\otimes\Leb,f\rangle}\bigr)
\end{equation}
This, along with \eqref{E:4.23}, completes the proof.
\end{proofsect}

\subsection{Distributional identity}
With the convergence \eqref{E:4.19} in hand, we are ready for the application of the coupling \eqref{E:4.5} which links every subsequential weak limit of random measures $\{\vartheta_{N\ge1}$ to the measures describing the thick points of the DGFF:

\begin{lemma}
Given~$f\colon \overline D\times\R_+\to\R$ with compact support, denote
\begin{equation}
f^{\ast\Leb}(x,\ell):=\int_\R\textd h\, f\bigl(x,\ell+\tfrac{h^2}2\bigr).
\end{equation}
Then every subsequential weak limit~$\vartheta$ of random measures $\{\vartheta_N\}_{N\ge1}$ satisfies
\begin{equation}
\label{E:4.31}
\bigl\langle\vartheta,f^{\ast\Leb}\bigr\rangle\,\, \laweq\,\,\int Z^D_{\sqrt\theta}(\textd x)\otimes\texte^{\alpha\sqrt\theta h}\textd h\, f\bigl(x,\tfrac{h^2}2\bigr)
\end{equation}
simultaneously for all~$f$ as above.
\end{lemma}

\begin{proofsect}{Proof}
Consider the coupling of~$L_t$, $h$ and~$\tilde h$ such that \eqref{E:4.5} holds. Denote
\begin{equation}
K_N:=\frac{N^2}{\sqrt{\log N}}\texte^{-\frac{(\sqrt{2t_N})^2}{2g\log N}}
\end{equation}
and observe that
\begin{equation}
\frac{\sqrt{\log N}}{\wh K_N}=\frac1{K_N}
\end{equation}
Given~$f$ as above, abbreviate
\begin{equation}
f^\ext(x,\ell,h):=f\bigl(x,\ell+\tfrac{h^2}2\bigr)
\end{equation}
The coupling then gives
\begin{equation}
\label{E:4.35}
\begin{aligned}
\frac1{K_N}\sum_{x\in D_N} f\bigl(x/N,\tfrac12 (\tilde h_x+\sqrt{2t_N})^2\bigr)
&=\frac{\sqrt{\log N}}{\wh K_N}\sum_{x\in D_N} f\bigl(x/N,L_t(x)+\tfrac12 h_x^2\bigr)
\\
&=\bigl\langle\vartheta_N^\ext,f^\ext\bigr\rangle
\end{aligned}
\end{equation}
Lemma~\ref{lemma-4.5} then tells us that, along the subsequence $\{N_k\}_{k\ge1}$ that takes $\vartheta_{N_k}$ to~$\vartheta$, the right-hand side tends weakly to
\begin{equation}
\langle\vartheta\otimes\Leb,f^\ext\rangle = \int\vartheta(\textd x\textd\ell)\,\textd h \,f\bigl(x,\ell+\tfrac{h^2}2\bigr) = \bigl\langle\vartheta,f^{\ast\Leb}\bigr\rangle
\end{equation}
On the other hand, since $\sqrt{2t_N}\sim 2\sqrt g\sqrt\theta\log N$, Theorem~\ref{thm-1.5} tells us that the left-hand side of \eqref{E:4.35} tends weakly to
\begin{equation}
\int Z^D_{\sqrt\theta}(\textd x)\otimes\texte^{\alpha\sqrt\theta h}\textd h\, f\bigl(x,\tfrac{h^2}2\bigr)
\end{equation}
which gives the desired claim.
\end{proofsect}

We now claim that this gives:

\begin{corollary}
Suppose $\mu$ is a deterministic measure on~$\R_+$ with the Laplace transform
\begin{equation}
\label{E:4.38}
\int_{\R+}\mu(\textd\ell)\texte^{-s\ell} = \exp\Bigl\{\frac{\alpha^2\theta}{2s}\Bigr\},\quad s>0
\end{equation} 
Then every subsequential weak limit~$\vartheta$ of random measures $\{\vartheta_N\}_{N\ge1}$ takes the form
\begin{equation}
\vartheta(\textd x\textd\ell) = Z^D_{\sqrt\theta}(\textd x)\otimes \mu(\textd\ell)
\end{equation}
\end{corollary}

\begin{proofsect}{Proof}
Given an open set~$A\subseteq\R^2$, abbreviate $\zeta_A(B):=\vartheta(A\times B)$. Fix~$s>0$ and denote $g_s(\ell):=\texte^{-s\ell}$. Take a sequence~$\{f_n\}_{n\ge1}$ of continuous compactly-supported functions that increase to $f:=1_A\otimes g$ and note that $f_n^{\ast\Leb}$ then increases to
\begin{equation}
f^{\ast\Leb}(x,\ell)=1_A(x)\texte^{-s\ell}\sqrt{\frac{2\pi}s}
\end{equation}
by the Monotone Convergence Theorem.
Noting that the equality in law in \eqref{E:4.31} was actually proved as an almost sure equality in a suitable coupling, applying the identity along the above sequence with the help of the Monotone Convergence Theorem shows
\begin{equation}
\sqrt{\frac{2\pi}s}
\int\zeta_A(\textd\ell)\texte^{-s\ell} = Z^D_{\sqrt\theta}(A)\int \textd h\,\texte^{\alpha\sqrt\theta h}\texte^{-s\frac12 h^2}\quad\text{a.s.}
\end{equation}
The Gaussian integral on the right equals $\texte^{\frac{\alpha^2\theta}{2s}}\sqrt{\frac{2\pi}s}$ which tells us that the measure
\begin{equation}
\nu_s(A):=\exp\Bigl\{-\frac{\alpha^2\theta}{2s}\Bigr\}\int_{A\times\R_+}\vartheta(\textd x\textd\ell)\texte^{-s\ell}
\end{equation}
equals to~$Z^D_{\sqrt\theta}(A)$ a.s. for all open~$A\subseteq\R^2$, regardless of~$s>0$. To overcome the fact that the implicit null set may depend on~$A$ and~$s$, note that the Borel sets in~$\R^2$ are generated by countably many open sets. A standard extension argument then gives that, on a set of full probability, $\nu_s = Z^D_{\sqrt\theta}$ for all rational~$s>0$. This means that, on a set of full measure,
\begin{equation}
\int_{A\times\R_+}\vartheta(\textd x\textd\ell)\texte^{-s\ell} = Z^D_{\sqrt\theta}(A)\int\mu(\textd\ell)\texte^{-s\ell}
\end{equation}
for all Borel~$A\subseteq \R^2$ and all~$s>0$, where continuity in~$s$ was used to get this for irrational~$s$. Since the Laplace transform determines the measure, the claim follows.
\end{proofsect}

We have basically proved:

\begin{theorem}
\label{thm-4.8}
For all $\theta\in(0,1)$ and any $\{t_N\}_{N\ge1}$ with~$t_N\sim 2g\theta(\log N)^2$,
\begin{equation}
\vartheta_N\,\,\underset{N\to\infty}\lawarrow\,\,Z^D_{\sqrt\theta}(\textd x)\otimes \mu(\textd\ell)
\end{equation}
where~$\mu$ is the measure
\begin{equation}
\mu(\textd\ell):=\delta_0(\textd\ell)+\biggl(\,\sum_{n\ge0}\frac1{n!(n+1)!}\Bigl(\frac{\alpha^2\theta}2\Bigr)^{n+1}\ell^n\biggr)1_{[0,\infty)}(\ell)\,\textd\ell
\end{equation}
\end{theorem}

\begin{proofsect}{Proof}
We just need to check that~$\mu$ has the Laplace transform \eqref{E:4.38} which is a straightforward calculation.
\end{proofsect}

\subsection{Proof of Theorem~\ref{thm-1.7}}
Intuitively, the measure~$\mu$ above captures the ``distribution'' of $O(1)$-values of the local time so the setting of Theorem~\ref{thm-1.7} should correspond to taking just the atom at~$0$ from~$\mu$. However, to make this precise we need to check that no part of the atom has come from infinitesimal values accumulating to zero in the limit. This is done in:

\begin{lemma}
\label{lemma-4.9}
For any~$\delta>0$ there exists~$c>0$ such that
\begin{equation}
\frac1{\wh K_N} \,E^\varrho\biggl(\,\,\sum_{\begin{subarray}{c}
x\in D_N \\d_\infty(x,D_N^\cc)>\delta N
\end{subarray}}
1_{\{0<L_{t_N}(x)\le\epsilon\}}\biggr)\le c\epsilon
\end{equation}
holds for all~$N\ge1$ and all~$\epsilon\in(0,1)$.
\end{lemma}

\begin{proofsect}{Proof}
Invoking one more time the representation \eqref{E:4.10} we have
\begin{equation}
P^\varrho\bigl(0<L_{t}(x)\le\epsilon\bigr)\le P\biggl(N_t\ge1\,\wedge\,\forall k=1,\dots, N_t\colon \sum_{k=1}^{Z_K}T_{k,j}\le\epsilon\pi(x)\biggr)
\end{equation}
Lemmas~\ref{lemma-4.1} and~\ref{lemma-4.3} then give
\begin{equation}
P^\varrho\bigl(0<L_{t}(x)\le\epsilon\bigr)\le \texte^{-\frac{t}{G(x,x)}\exp\{-\frac{\epsilon}{G(x,x)}\}}-\texte^{-\frac{t}{G(x,x)}}
\le \texte^{-\frac{t}{G(x,x)}}\bigl(\texte^{\epsilon\frac{t}{G(x,x)^2}}-1\bigr)
\end{equation}
Now set~$V:=D_N$ and $t_N=O((\log N)^2)$. For any~$\delta>0$ small, once~$x$ at least~$\delta N$ from the boundary, we have $t_N/G^{D_N}(x,x)^2=O(1)$ and so
\begin{equation}
P^\varrho\bigl(0<L_{t}(x)\le\epsilon\bigr)\le c\frac{\wh K_N}{N^2}\epsilon
\end{equation}
The claim follows by summing this over $x\in D_N$ with $d_\infty(x,D_N^\cc)>\delta N$.
\end{proofsect}

We now finally give:

\begin{proofsect}{Proof of Theorem~\ref{thm-1.7}}
Let~$g\colon D\to\R$ be continuous with compact support and, for each~$n\ge1$, set $f_n(x,\ell):=g(x)(1-n\ell)_+$. Denote
\begin{equation}
\kappa_N:=\frac1{\wh K_N}\sum_{x\in D_N}1_{\{L_{t_N}(x)=0\}}\delta_{x/N}
\end{equation}
Then 
\begin{equation}
\bigl|\langle \vartheta_N,f_n\rangle-\langle\kappa_N,g\rangle\bigr|\le\Vert g\Vert
\,\frac1{\wh K_N}\!\!
\sum_{\begin{subarray}{c}
x\in D_N \\ x/N\in\supp(g)
\end{subarray}}
1_{\{0<L_{t_N}(x)\le1/n\}}
\end{equation}
Since $d_\infty(\supp(g),D^\cc)>\delta$, Lemma~\ref{lemma-4.9} tells us that
\begin{equation}
\lim_{n\to\infty}\limsup_{N\to\infty}\, E^\varrho\bigl|\langle \vartheta_N,f_n\rangle-\langle\kappa_N,g\rangle\bigr|=0
\end{equation}
Invoking Theorem~\ref{thm-4.8}, this gives
\begin{equation}
\langle\kappa_N,g\rangle\,\,\underset{N\to\infty}\lawarrow\,\, \bigl\langle Z^D_{\sqrt\theta},g\bigr\rangle\,\lim_{n\to\infty}\int \mu(\textd\ell)f_n(\ell)
\end{equation}
The limit on the right equals~$1$ and so~$\kappa_N\,\lawarrow\, Z^D_{\sqrt\theta}$ as measures on~$D\times\R_+$. The convergence is extended to~$\overline D\times\R_+$ with the help of tightness proved in Corollary~\ref{cor-4.4} and the fact that~$\vartheta_N$ naturally dominates~$\kappa_D$.
\end{proofsect}

We finish by noting that the joint paper~\cite{AB} with Y.~Abe where the above results are proved contains results for other level sets; namely, the \textit{$\lambda$-thick points} of~$L_{t_N}$,
\begin{equation}
\bigl\{x\in D_N\colon L_{t_N}(x)\ge2g(\sqrt\theta+\lambda)^2(\log N)^2\bigr\}
\end{equation}
where~$\lambda\in(0,1)$, and~the \textit{$\lambda$-thin points} of~$L_{t_N}$,
\begin{equation}
\bigl\{x\in D_N\colon L_{t_N}(x)\ge2g(\sqrt\theta-\lambda)^2(\log N)^2\bigr\}
\end{equation}
where~$\lambda\in(0,\min\{1,\sqrt\theta\})$ and~$t_N\sim 2g\theta(\log N)^2$. The spatial distribution of these points is again described by the measure~$Z^D_\lambda$. All of these results demonstrate universality of the DGFF in these extremal problems.

\renewcommand{\sectionphrase}{Appendix}

\newpage
\section{What lies beyond?}
\noindent
In the lectures we concentrated on problems whose statement is easy to make and the proof can realistically be presented in the available time. In order to make these notes more useful, we will now comment on results that were omitted and point out a few interesting conjectures in this subject area. We will nonetheless remain focused and so our presentation should not be regarded as a comprehensive review.

\subsection{Maximum and extremal process of DGFF}
The control of the thick points of the DGFF naturally leads to the question of the limit law of the absolute maximum and the associated extremal process. For the maximum one gets the limit theorem
\begin{equation}
\label{E:5.1}
\BbbP\Bigl(\,\max_{x\in D_N} h^{D_N}_x\le m_N+u\Bigr)\,\,\underset{N\to\infty}\longrightarrow\,\, \E\bigl(\texte^{-\eusb Z\texte^{-\alpha u}}\bigr)
\end{equation}
where $\alpha:=2/\sqrt g$, the centering is done by the sequence
\begin{equation}
\label{E:5.2}
m_N:=2\sqrt g\,\log N-\frac34\sqrt g\,\log\log N
\end{equation}
and $\eusb Z$ is an a.s.-finite and positive random variable. The centered maximum thus tends to a \textit{randomly-shifted Gumbel law}, i.e.,
\begin{equation}
\max_{x\in D_N} h^{D_N}_x - m_N\,\,\,\underset{N\to\infty}\lawarrow\,\,\, G+\alpha^{-1}\log\eusb Z
\end{equation}
where~$G$ is a Gumbel random variable independent of~$\eusb Z$. 

The limit \eqref{E:5.1} was proved by M.~Bramson, J.~Ding and O.~Zeitouni~\cite{BDingZ} building upon an earlier proof of tightness by M.~Bramson and O.~Zeitouni~\cite{BZ}. Both of these involve comparisons with modified Branching Random Walk. A different argument is presented in~\cite[Lectures 8-10]{B-notes} that uses the Gibbs-Markov property, a novel inequality expanding on the Dekking-Host argument~\cite{Dekking-Host} and a concentric decomposition of the DGFF. 

A natural next question is: What is~$\eusb Z$? And, can it be characterized independently? This was  settled only in a sequence of papers~\cite{BL1,BL2,BL3} by O.~Louidor and the author of these notes. The interest there is on the \textit{extremal process} 
\begin{equation}
\eta_N:=\sum_{x\in D_N}\delta_{x/N}\otimes\delta_{h_x^{D_N}-m_N}
\end{equation}
which is similar to the processes associated with the $\lambda$-thick points, albeit without any normalization. For this process, the above papers showed that, given any compactly-supported continuous $f\colon\overline D\times\R\to\R_+$, we have
\begin{equation}
\label{E:5.5}
\E\bigl(\texte^{-\langle\eta_N,f\rangle}\bigr)
\,\underset{N\to\infty}\longrightarrow\,\E\Biggl(\exp\biggl\{-\int Z^D(\textd x)\otimes\texte^{-\alpha h}\textd h\otimes\DD(\textd\zeta)(1-\texte^{-\langle\zeta,f(x,h-\cdot)\rangle})\biggr\}\Biggr)
\end{equation}
Here the objects on the right are as follows:
\settowidth{\leftmargini}{(111)}
\begin{itemize}
\item $Z^D$ is an a.s.-finite random Borel measure on~$D$ that is non-atomic and charges non-empty every open subset of~$D$ almost surely, and
\item $\DD$ is a (non-random) probability measure on the set of point measures (i.e., $\N\cup\{+\infty\}$-valued Radon measures) on~$\R$ such that a.e.-sample~$\zeta$ from~$\DD$ is infinite, supported on $(-\infty,0]$ and has a single atom at zero.
\end{itemize}
The expectation on the right of \eqref{E:5.5} is over~$Z^D$.

\begin{figure}[t]
\refstepcounter{obrazek}
\centerline{\includegraphics[width=5.6in]{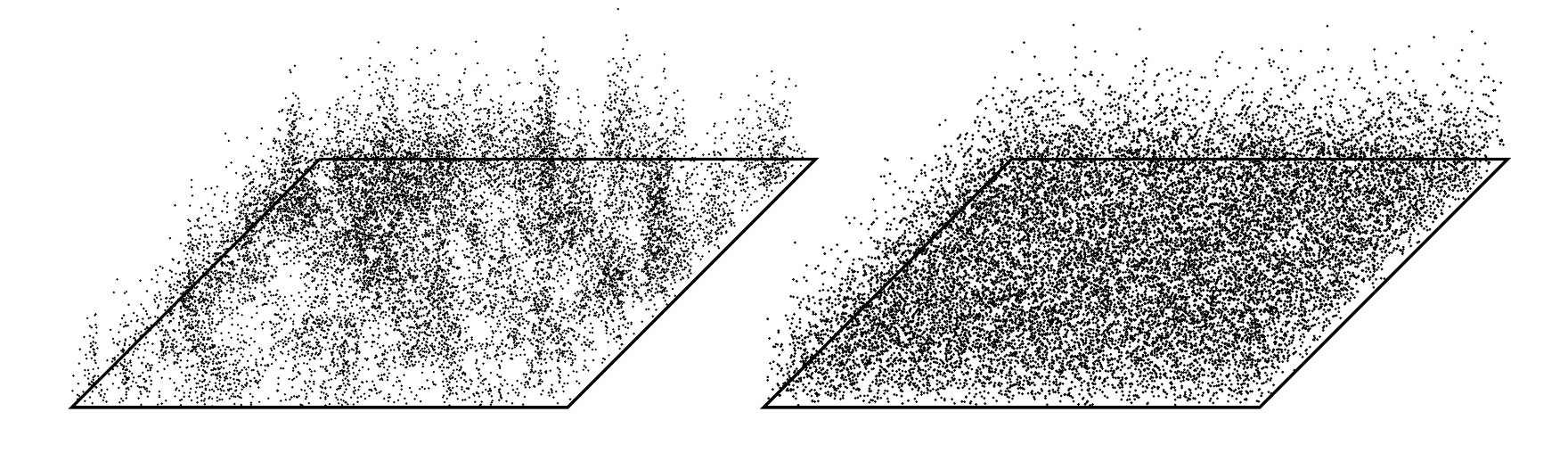}}
\begin{quote}
\label{fig10}
\fontsize{9}{5}\selectfont
{\bf Fig.~\theobrazek:\ }A comparison between the large values of the DGFF (left) and i.i.d.\ Gaussians (right) demonstrating the clustering of the DGFF values.
\normalsize
\medskip
\end{quote}
\end{figure}

\smallskip
As the above may be difficult to parse, here is another way to state the limit result: Let $\{(x_i,h_i)\}_{i\ge1}$ enumerate points in a sample from the Poisson point process on~$D\times\R$ with random intensity measure
\begin{equation}
Z^D(\textd x)\otimes\texte^{-\alpha h}\textd h
\end{equation}
and, given i.i.d.\ samples $\{\zeta^{(i)}\}_{i\ge1}$ from~$\DD$, let $\{d_j^{(i)}\}_{j\ge1}$ enumerate the points in~$\zeta^{(i)}$. Then
\begin{equation}
\label{E:5.7}
\eta_N\,\,\,\underset{n\to\infty}\lawarrow\,\,\,\sum_{i\ge1}\sum_{j\ge1}\delta_{x_i}\otimes\delta_{h_i+d_j^{(i)}}
\end{equation}
Each point $(x_i,h_i)$ thus carries with it a whole cluster of points at ``heights'' $h_i+d^{(i)}_j$ with the highest point in the cluster at height~$h_i$ representing a local maximum of DGFF. (The clusters are sometimes referred to as ``decorations'' and the process on the right of \eqref{E:5.7} is called a \textit{decorated Poisson point process}.) See Fig.~\ref{fig10} for a visualization of clustering of extremal DGFF values.  


A simple approximation argument shows that \eqref{E:5.5} implies \eqref{E:5.1} with~$\eusb Z$ corresponding to the total mass of~$Z^D$ measure,
\begin{equation}
\eusb Z\,\laweq\, \alpha^{-1}Z^D(D)
\end{equation}
but we still owe the reader a characterization of~$Z^D$ itself. As it turns out, $Z^D$ is a critical case of the Liouville Quantum Gravity (LQG) corresponding, roughly, to the $\lambda=1$ boundary case of the family~$\{Z^D_\lambda\colon \lambda\in(0,1)\}$ discussed in Lecture~3. (The $\lambda=0$ boundary case is given by the Lebesgue measure.)

The critical nature of $Z^D$ stems from the fact that the construction along the lines \eqref{E:3.65} will yield a zero measure; one has to scale-up the approximations by a suitable~$n$-dependent factor  to obtain a non-trivial limit. This is related to the fact that~$Z^D(A)$ does not admit a first moment; indeed, for all~$A\subseteq\R^2$ open we instead have
\begin{equation}
\lim_{t\downarrow0}\,
\frac1{ \log(1/t)}
E\bigl(Z^D(A)\,\texte^{-t Z^D(D)}\bigr)\,\,=\bar c\int_A r^D(x)^2\textd x
\end{equation}
for a constant~$\bar c\in(0,\infty)$. As it turns out, this along with the Gibbs-Markov property
\begin{equation}
Z^D(\textd x)\laweq\,\texte^{\alpha\Phi^{D,\wt D}(x)}Z^{\wt D}(\textd x),\quad \Phi^{D,\wt D}\independent Z^{\wt D}
\end{equation}
whenever~$\wt D\subseteq D$ with~$\Leb(D\smallsetminus\wt D)=0$
 nails the measure uniquely; i.e., the only freedom one has is the choice of constant~$\bar c$. See~\cite[Chapter~10]{B-notes} for more details.

The critical LQG measure can be constructed for a fairly large class of  logarithmically correlated Gaussian fields. The constructions go back to B.~Duplantier, R.~Rhodes, S.~Sheffield and V.~Vargas~\cite{DRSV1,DRSV2} with uniqueness settled in full generality by J.~Junnila and E.~Saksman~\cite{JS} and E.~Powell~\cite{Powell}. As shown in a recent preprint~\cite{BGL} by S.~Gufler, O.~Louidor and the author,~$Z^D$ admits a carrier of vanishing Hausdorff dimension almost surely. Supercritical (i.e., $\lambda>1$) variants of the LQG measure can be constructed as well, but they are purely atomic. See~R.~Rhodes and V.~Vargas~\cite{RV-review} or~\cite[Theorem~2.6]{BL3}.

\subsection{More on local time thick points}
As explained in Lecture~2, the local time of two dimensional random walk shares many features of the DGFF (in particular, it is also a logarithmically correlated field) and so it presents us with a natural playground to attempt proofs of analogous conclusions. However, even the situation of the thick points has not yet been completely resolved. Indeed, Theorem~\ref{thm-1.7} worked under a return mechanism based on the boundary vertex~$\varrho$, and it parametrized time by the local time at~$\varrho$.

The conversion to natural time parametrization has been studied by Y.~Abe, S.~Lee and the author in~\cite{ABL}. To state the result, fix an admissible~$D\subseteq\R^2$ and let~$\{D_N\}_{N\ge1}$ be admissible approximations of~$D$. We first need an upgrade of Theorem~\ref{thm-1.5} that allows for conditioning on zero average of~DGFF on~$D_N$:

\begin{theorem}
\label{thm-5.1}
There exists a family of random Borel measures~$\{Z^{D,0}_\lambda\colon\lambda\in(0,1)\}$ such that for any $\{a_N\}_{N\ge1}$ with $a_N\sim 2\sqrt g\lambda\log N$ for some $\lambda\in(0,1)$ and~$K_N$ defined from~$a_N$ as in \eqref{E:1.26},
\begin{equation}
\biggl(\frac1{K_N}\sum_{x\in D_N}\delta_{x/N}\otimes\delta_{h_x^{D_N}-a_N}\,\bigg|\,\sum_{x\in D_N}h_x^{D_N}=0\biggr)\,\,\,\underset{N\to\infty}\lawarrow\,\,\,
Z^{D,0}_\lambda(\textd x)\otimes\texte^{-\alpha \lambda h}\textd h
\end{equation}
Moreover,
if $Y$ is a normal random variable with mean zero and variance
\begin{equation}
\Var(Y):=\int_{D\times D}\textd x\textd y\,\wh G^D(x,y)
\end{equation}
then
\begin{equation}
\label{E:2.15new}
Y\independent Z^{D,0}_\lambda\quad\Rightarrow\quad
\texte^{\lambda\alpha\frakd(x)Y}
\,Z^{D,0}_\lambda(\textd x)\,\laweq\,Z^D_\lambda(\textd x)
\end{equation}
where
\begin{equation}
\frakd(x):=\Leb(D)\frac{\int_D\textd y\,\wh G^D(x,y)}{\int_{D\times D}\textd z\,\textd y\,
\wh G^D(z,y)}
\end{equation}
The law of $Z^{D,0}_\lambda$ is determined uniquely by \eqref{E:2.15new}.
\end{theorem}

\begin{figure}[t]
\refstepcounter{obrazek}
\centerline{\includegraphics[width=3.3in]{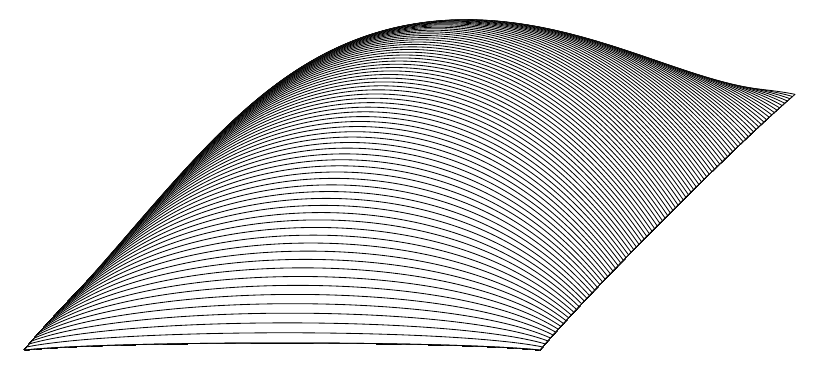}}
\begin{quote}
\label{fig11}
\fontsize{9}{5}\selectfont
{\bf Fig.~\theobrazek:\ }A plot of function~$\frakd$ on~$D:=(0,1)\times(0,1)$ obtained by solving the Poisson equation with constant charge density $-\Delta\frakd = \Leb(D)/\Var(Y)$.
\normalsize
\medskip
\end{quote}
\end{figure}

Moving to the local time let $\deg(\overline D_N)$ denote the sum of all degrees of the vertices in~$D_N\cup\{\varrho\}$ (or, alternatively, twice the number of all edges). The natural time parametri\-zation of the local time is then achieved by taking
\begin{equation}
\label{E:5.15}
\wt L_t(x):=\ell_{t\deg(\overline D_N)}(x)
\end{equation}
The result for the avoided points then reads:

\begin{theorem}
\label{thm-5.2}
Suppose that $\{t_N\}_{N\ge1}$ obeys $t_N\sim 2g\theta\log N$ for some~$\theta\in(0,1)$ and define~$\wh K_N$ from~$t_N$ as in \eqref{E:1.33}. Then
\begin{equation}
\frac1{\wh K_N}\sum_{x\in D_N} 1_{\{\wt L_{t_N}(x)=0\}}\delta_{x/N}\,\,\,\underset{N\to\infty}\lawarrow\,\,\,\,\texte^{\alpha\sqrt\theta(\frakd(x)-1)Y}Z^{D,0}_{\sqrt\theta}
\end{equation}
where $Y$ and~$Z^{D,0}_{\sqrt\theta}$ are independent and with laws as specified in Theorem~\ref{thm-5.1}.
\end{theorem}

The upshot is that, in the natural time parametrization, the spatial distribution of the avoided points at $\theta$-multiple of the cover time is described by a measure somewhere ``in-between'' the limit spatial distribution of the $\sqrt\theta$-thick points of~$h^{D_N}$ and of~$h^{D_N}$ conditioned on vanishing $\sum_{x\in D_N}h^{D_N}_x$. The function~$\frakd$ shows that the effect of the conditioning distributes inhomogeneously throughout~$D$; see Fig.~\ref{fig11}.

Both theorems are deduced from the earlier conclusions by way of manipulations that, besides the coupling in Theorem~\ref{thm-2.7}, require solving several non-trivial distributional identities. See~\cite{ABL} for details.

The above results are still deficient due to their reliance on the boundary vertex (which is what makes Theorem~\ref{thm-2.7} available). A natural setting for the study of the local time are domains with free boundary conditions and the lattice torus. These cases are currently subject of various research attempts.

\subsection{Frequent points and cover time scaling}
There are two natural extremal problems associated with the local time that have not yet been resolved. The first one concerns the \textit{frequent points} which are those where the local time (associated with a path of the random walk) is maximal. The second problem concerns the precise scaling, and a limit law, of the \textit{cover time}. 

\begin{figure}[t]
\refstepcounter{obrazek}
\centerline{\includegraphics[width=4.9in]{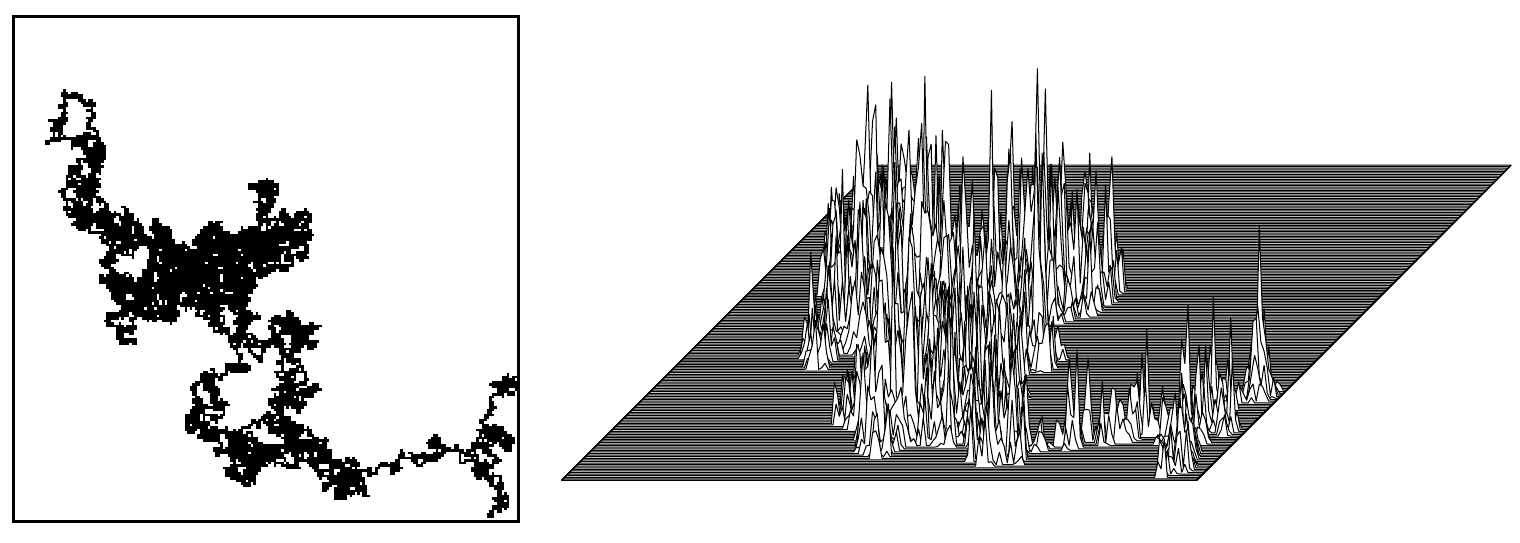}}
\begin{quote}
\label{fig12}
\fontsize{9}{5}\selectfont
{\bf Fig.~\theobrazek:\ }A plot of the range (left) and the local time configuration (right) of the simple random walk started at the origin of a lattice square and run until the first exit.
\normalsize
\medskip
\end{quote}
\end{figure}

The study of frequent points was initiated in an influential paper by P.~Erd\H os and S.J.~Taylor~\cite{ET60} who identified the two-dimensional case as the one most interesting. Even  the leading order for the latter case took full 40 years to be resolved by A.~Dembo, Y.~Peres, J.~Rosen and O.~Zeitouni~\cite{DPRZ01}. Casting the problem in the context of our scaled-up versions~$D_N$ (see Fig.~\ref{fig12}) of an admissible domain~$D$, which we assume contains the origin in~$\R^2$, with the complement of~$D_N$ collapsed to the boundary vertex~$\varrho$, the question is:
\begin{equation}
\text{What is the asymptotic limit law of $\max_{x\in D_N}\ell_{\tau_\varrho}(x)$ under~$P^0$?}
\end{equation}
The leading order asymptotic is known to be
\begin{equation}
\frac1{(\log N)^2}\max_{x\in D_N}\ell_{\tau_\varrho}(x)\,\,\underset{N\to\infty}{\overset{P}\longrightarrow}\,\,\frac 4\pi
\end{equation}
The task is to find all the subleading deterministic orders up to the order at which a non-degenerate weak limit is possible. The following is expected to hold (see, e.g., A.~Jego~\cite[Conjecture~1]{J20}):

\begin{conjecture}
\label{conj-5.3}
There exists $\bar\alpha\in(0,\infty)$ and, for each admissible domain~$D$ containing the origin, there exists an a.s.-finite and positive random variable~$\eusb Z^D$ such that
\begin{equation}
P^0\biggl(\,\max_{x\in D_N}\ell_{\tau_\varrho}(x)\le\frac2{\sqrt\pi}\log N-\frac1{\sqrt\pi}\log\log N+u\biggr)\,\,\underset{N\to\infty}\longrightarrow\,\,E\bigl(\texte^{-\eusb Z^D\texte^{-\bar\alpha u}}\bigr)
\end{equation}
holds for each~$u\in\R$.
\end{conjecture}

\noindent
A similar conjecture can also be made about the associated extremal process which, just as for the DGFF, is expected to have the structure of the decorated Poisson point process with a random intensity measure; the random variable~$\eusb Z^D$ is then the total mass of this measure times~$\bar\alpha^{-1}$. In~\cite{J20}, A.~Jego actually constructed a candidate for the random intensity measure called the \textit{Brownian multiplicative chaos} drawing on measures that arise in his study of the local-time thick points in~\cite{J18}.

Progress on proving Conjecture~\ref{conj-5.3} has taken place mainly in other contexts. Indeed, J.~Rosen~\cite{Rosen-new} proved tightness of the maximum around the above scale sequence, albeit for the Brownian motion on a sphere in~$\R^2$ and the local time interpreted via the average occupation time measure on balls of radius~$1/N$. O.~Louidor and the author~\cite{BL5}, based on earlier work of~Y.~Abe~\cite{A18}, proved the conjecture for the random walk on a $b$-ary rooted tree of depth~$n$, with the walk started at a leaf-vertex and stopped when the root it hit. The latter context is more amenable to analysis because the local time parametrized by the time spent at the root (to which the stated conclusion can be converted) has a tree-indexed Markovian structure.

\smallskip
The situation with the cover time is roughly similar. The goal is to show that a suitably centered and scaled cover time of a domain of linear size~$N$ in~$\Z^2$ or the lattice torus $(\Z/N\Z)^2$ tends to a non-degenerate limit. (The leading order and the case of other dimensions was discussed in and after \eqref{E:1.8}.) Focusing on the easier case of simple random walk in admissible approximations~$D_N$ of admissible domains, with returns done via the boundary vertex~$\varrho$, here one expects the following:

\begin{conjecture}
\label{conj-5.4}
For each admissible domain~$D$, there exists an almost-surely finite and positive random variable~$\overline{\eusb Z}^D$ such that
\begin{equation}
\label{E:5.20}
P^\varrho\Biggl(\,\sqrt{\frac{\tau_{\text{\rm cov}}}{\deg(\overline D_N)}}\le \frac1{\sqrt\pi}\log N -\frac1{4\sqrt\pi}\log\log N+u\Biggr)
\,\,\underset{N\to\infty}\longrightarrow\,\,E\bigl(\texte^{-\overline{\eusb Z}^D\texte^{-\sqrt{4\pi}\, u}}\bigr)
\end{equation}
holds for all~$u\in\R$. Here $\text{\rm deg}(\overline D_N)$ is as in \eqref{E:5.15}.
\end{conjecture}

This conjecture can largely be explained by the beautiful argument underlying the recent work of O.~Louidor and S.~Saglietti~\cite{LS-new} who used it to prove that the CDFs on the left of \eqref{E:5.20} form a tight family. 

The bulk of the argument relies on the time parametrization at the boundary vertex. The idea is to split the (expected) time scale of the cover time in two phases. In the first phase, we will first run the random walk for time
\begin{equation}
t_A:=\frac12 m_N^2
\end{equation}
 where~$m_N$ is as in \eqref{E:5.2}. The reason for this choice is the fact that, in the coupling from Theorem~\ref{thm-2.7}, for each~$a>0$ we then have
\begin{equation}
L_{t_A}(x)\le a^2\,\wedge\, |h_x|\le \sqrt 2\,a\,\,\Rightarrow\,\, |\tilde h_x+m_N|\le 2a
\end{equation}
Note that \eqref{E:5.7} tells us that  $\{x\in D_n \colon |\tilde h_x+m_N|\le 2a\}$ is typically composed of a finite number of  components separated by distances of order~$N$, and the total number of such components is proportional to the total mass of~$Z^D$-measure.  Invoking the thinning argument underlying our derivations in Lecture~4, it is thus reasonable to expect that $\{x\in D_N\colon L_{t_A}(x)=0\}$ is composed of $c Z^D(D)\sqrt{\log N}$ (mostly) small islands well separated (typically by distance $N^{1-o(1)}$) from each other.

After the first phase is done we ``restart'' the random walk and run it for time $t_B+(\sqrt{2g}\log N)u$, where
\begin{equation}
t_B:=\frac12(g\log N)\log\log N
\end{equation}
In this phase, the aim is to ensure that each of the components of $\{x\in D_N\colon L_{t_A}(x)=0\}$ is visited by at least one excursion from~$\varrho$. (We are assuming that once the component is hit, it will be swept out completely.) By the calculation in the proof of Lemma~\ref{lemma-4.1}, the number of excursions that hit vertex~$x$ deep inside~$D_N$ is Poisson with parameter
\begin{equation}
\frac{t_B+(\sqrt{2g}\log N)u}{G(x,x)}
=\frac12\log\log N+\sqrt{2/g}\, u+o(1)
\end{equation}
and so the probability that a given component is hit equals $\texte^{-\frac12\log\log N+\sqrt{2/g}\, u+o(1)}$. Assuming that each excursion will typically hit only one component, all $c Z^D(D)\sqrt{\log N}$ components of $\{x\in D_N\colon L_{t_A}(x)=0\}$ are hit by time $t_B+(\sqrt{2g}\log N)u$ with probability 
\begin{equation}
\Bigl(1-\texte^{-\frac12\log\log N-\sqrt{2/g}\, u+o(1)}\Bigr)^{c Z^D(D)\sqrt{\log N}}
=\texte^{-c Z^D(D)\texte^{-\sqrt{2/g}\, u}+o(1)}
\end{equation}
Since $g=(2\pi)^{-1}$, a calculation now shows
\begin{equation}
t_A+t_B+(\sqrt{2g}\log N)u = \Bigl(\frac1{\sqrt\pi}\log N-\frac1{4\sqrt\pi}\log\log N+u+o(1)\Bigr)^2
\end{equation}
Denoting the cover time in time parametrization at the boundary vertex as
\begin{equation}
\wt\tau_{\text{\rm cov}}:=\inf\Bigl\{t\ge0\colon \min_{x\in D_N\cup\{\varrho\}}L_t(x)>0\Bigr\}
\end{equation}
we thus conclude
\begin{equation}
P^\varrho\Biggl(\,\sqrt{\wt\tau_{\text{\rm cov}}}\le \frac1{\sqrt\pi}\log N -\frac1{4\sqrt\pi}\log\log N+u\Biggr)
\,\,\underset{N\to\infty}\longrightarrow\,\,E\bigl(\texte^{-cZ^D(D)\texte^{-\sqrt{4\pi}\, u}}\bigr)
\end{equation}
where we also observed that $\sqrt{2/g}=\sqrt{4\pi}$.

This is unfortunately not the end of the story because we still need to make a conversion to the actual time. The overall rescaling by~$\deg(\overline D_N)$ is as in \eqref{E:5.15}, but there is also a correction due to the fluctuations that made the transition from Theorem~\ref{thm-1.7} to Theorem~\ref{thm-5.2} worthy of a paper~\cite{ABL}. We thus expect~$\overline{\eusb Z}^D$ to be a constant multiple of the total mass of $\theta=1$ version of the measure on the right of \eqref{E:5.15}, rather than~$Z^D$ itself.

\smallskip
A corresponding version of Conjecture~\ref{conj-5.4} was proved along the above lines for the random walk on a $b$-ary rooted tree in the paper by A.~Cortinez, O.~Louidor and S.~Saglietti~\cite{CLS}. Tightness of the (suitably defined) cover time for Brownian motion on a two-dimensional sphere was shown by D.~Belius, J.~Rosen and O.~Zeitouni~\cite{BRZ}.

\newpage

\centerline{\bf\large Acknowledgements}
\vglue2mm\noindent
It is a pleasure to thank the organizers of the school, specifically, \'Agnes Backhausz, G\'abor Pete, Bal\'azs R\'ath and B\'alint T\'oth, for the opportunity to share the results and views of this subject area. The author's travel to the school and some of the reported work have been  supported by the NSF awards DMS-1954343 and DMS-2348113.

\bigskip

\renewcommand\refname{\bf References}

\refstepcounter{section}
\addcontentsline{toc}{section}
{{\tocsection {}{\thesection}
{\!\!\!References\dotfill}}{}}

\bibliographystyle{abbrv}


\end{document}